\documentclass[10pt,letterpaper]{article}
\usepackage[top=0.85in,left=2.75in,footskip=0.75in]{geometry}

\usepackage{amsmath,amssymb}

\usepackage{changepage}

\usepackage[utf8x]{inputenc}

\usepackage{textcomp,marvosym}

\usepackage{cite}

\usepackage{nameref,hyperref}

\usepackage[right]{lineno}

\usepackage{microtype}
\DisableLigatures[f]{encoding = *, family = * }

\usepackage[table]{xcolor}

\usepackage{array}

\newcolumntype{+}{!{\vrule width 2pt}}

\newlength\savedwidth

\raggedright
\usepackage[aboveskip=1pt,labelfont=bf,labelsep=period,justification=raggedright,singlelinecheck=off]{caption}

\makeatletter
\renewcommand{\@biblabel}[1]{\quad#1.}
\makeatother

\usepackage{lastpage,fancyhdr,graphicx}
\usepackage{epstopdf}
\fancyheadoffset[L]{2.25in}
\fancyfootoffset[L]{2.25in}
\usepackage{amsthm}

\newcommand{\ee}{\epsilon}
\newcommand{\U}{\mathcal{U}}
\newcommand{\A}{\mathcal{A}}
\newcommand{\M}{\mathcal{M}}
\newcommand{\RR}{\mathbb{R}}
\newcommand{\DMM}{\left[\partial \mathcal{M}\right]_-}
\newcommand{\DMO}{\left[\partial \mathcal{M}\right]_0}

\newcommand{\Acomp}{\mathcal{A}^{\textsf{C}}}
\newcommand{\DAM}{\left[\partial \mathcal{A}\right]_-}
\newcommand{\DAO}{\left[\partial \mathcal{A}\right]_0}
\newcommand{\ie}{\emph{i.e.}}
\newcommand{\eg}{\emph{e.g.}}
\newtheorem{theorem}{Theorem}[section]
\newtheorem{remark}{Remark}[section]
\newtheorem{definition}{Definition}[section]
\newtheorem{proposition}{Proposition}[section]
\newtheorem{problem}{Problem}

\begin{document}
\vspace*{0.2in}

\begin{flushleft}
{\Large
\textbf\newline{Epidemic Management with Admissible and Robust Invariant Sets} 
}
\newline
\\
Willem Esterhuizen\textsuperscript{1 *},
Jean Lévine\textsuperscript{2},
Stefan Streif\textsuperscript{1},
\\
\bigskip
\textbf{1} Technische Universität Chemnitz, Automatic Control and System Dynamics Laboratory, Germany
\\
\textbf{2} Centre Automatique et Systèmes (CAS), Maths \& Systems Dept., MINES ParisTech, PSL University, France
\\
\bigskip

%
%

%
%
%

* willem.esterhuizen@etit.tu-chemnitz.de

\end{flushleft}
\section*{Abstract}
We present a detailed set-based analysis of the well-known SIR and SEIR epidemic models subjected to hard caps on the proportion of infective individuals, and bounds on the allowable intervention strategies, such as social distancing, quarantining and vaccination. We describe the admissible and maximal robust positively invariant (MRPI) sets of these two models via the theory of barriers. We show how the sets may be used in the management of epidemics, for both perfect and imperfect/uncertain models, detailing how intervention strategies may be specified such that the hard infection cap is never breached, regardless of the basic reproduction number. The results are clarified with detailed examples.



\section{Introduction}

There is a large literature on the application of optimal control to epidemiology. Some of the earliest papers on the topic are \cite{sanders1971quantitative}, which investigates the optimal control of a disease of SIS type (susceptible-infective-susceptible), and \cite{hethcote1973optimal}, which considers an SIR (susceptible-infective-removed) model with vaccination to obtain optimal vaccination policies via dynamic programming. Other papers consider the optimal control of SIR models, \cite{HanD11,kruse2020optimal, miclo2020optimal, godara2021control}; of HIV \cite{kirschner1997optimal, culshaw2004optimal}; and of malaria \cite{agusto2017optimal}. This list is by no means exhaustive, and we point the reader to the survey paper, \cite{sharomi2017optimal}. Moreover, scores of papers have recently appeared involving the modelling and control of COVID-19, often via model-predictive control (MPC), see for example \cite{ARONNA2021100437, bonnans2020optimal, grundel2020much, perkins2020optimal, kohler2020robust} and the thorough discussion of \cite{Lobry2021}. 

In this paper we build on the recent work concerning the application of \emph{set-based methods} to epidemiology. We present a set-based analysis of two well-known \emph{compartmental epidemic models}: the SIR and SEIR models, see for example \cite{Brauer2019, Het2000}, subjected to constraints on their inputs (social distancing, quarantine rate and/or vaccination rate) and state (a hard cap on the proportion of infectives), utilising the \emph{theory of barriers}, \cite{DonL13,EstAS20}. These models often serve as good initial candidates to model new diseases, later being elaborated on to be more accurate.

It has only recently been shown that set-based methods may be used to maintain hard infection caps in epidemic models. To our knowledge, the paper \cite{LarS16} introduced the idea, describing the so-called \emph{viability kernel}\footnote{the set of all states for which there exists an input such that the state and input constraints are satisfied for all future time}, \cite{aubin2009viability}, of a two-dimensional model of a vector-borne disease. The work was then extended in \cite{SalL19} where the effects of modelling uncertainties of this system were investigated via the \emph{robust} viability kernel. The recent paper \cite{RASHKOV2021124958} finds the viability kernel for an SIR model of a vector bourne disease by approximating the value function of an appropriate optimal control problem by solving the related Hamilton-Jacobi-Bellmann (HJB) equation.

In our previous research, \cite{EstALetal20}, we showed that the theory of barriers may be utilised to describe the viability kernel (also known as the \emph{admissible set}) as well as the \emph{maximal robust positively invariant set} (MRPI)\footnote{the set of states in which the state and input constraints are satisfied for all time regardless of the input} of the malaria model considered in \cite{LarS16}. The theory of barriers describes the non-trivial part of these sets' boundaries, made of integral curves of the system that satisfy a minimum-/maximum-like principle, allowing them to easily be constructed. This contrasts with other approaches that estimate the viability kernel with algorithms that iteratively compute reachable sets, \cite[Ch. 5]{Blanchini_2015}, with interval analysis, \cite{Monnet_2016}, through the solution of HJB equations, \cite{Mitchell_2005}, or via polynomial optimisation, \cite{korda2014}, all of which suffer from the curse of dimensionality. Moreover, optimal control problems with hard state constraints (such as a hard cap on the number of infectives) are notoriously difficult to solve due to the presence of ``active arcs'', ``jump conditions'' and adjoint dynamics that are often difficult to specify. The construction of the sets via the theory of barriers does not suffer from these difficulties because the state constraints do not directly appear in their construction.

Other research that contributes to epidemiology with set-based ideas includes the paper \cite{barrios2018sustainable}, which describes the \emph{set of sustainable thresholds} (in some sense the dual of the viability kernel) for the class of so-called \emph{cooperative} epidemic models. An attractive aspect of this approach is that competing costs (for example, number of infective individuals versus the disease containment costs) may be analysed via Pareto efficiency. The work in \cite{kassara2009unified,boujallal2021novel} is closely related to set-based methods, but differs in that they do not describe or compute sets. Rather, these works combine Lyapunov-like functions with ideas from Viability theory (the so-called ``regulation map'') to specify continuous selections of set-valued maps that describe the elimination of a disease.

The set-based approach to the management of epidemics, as presented in \cite{LarS16,SalL19,EstALetal20}, argues that intervention strategies (for example, the level of fumigation in the control of mosquito-borne diseases) should be chosen based on the location of the state with respect to the computed sets: 
If the state is located in the admissible set, then it is possible to maintain the infection cap for all time with a suitably chosen input (intervention strategy). If the state is located in the MRPI, then the cap will never be breached and any input may be used, such as, \eg, saving resources or minimising economic damage.

In the current paper we analyse these sets for the constrained SIR and SEIR models under two cases: a \emph{perfect model case} and an \emph{imperfect model case}. For the perfect model case the modelling parameters are assumed to be perfectly known and the intervention strategy may be chosen depending on the location of the state, as was described earlier. For the imperfect model case we assume that there is a \emph{pre-designed intervention strategy}, possibly a feedback, that has been designed and implemented on the system. For this case we only describe the MRPI, which corresponds to states from which the infection cap can always be maintained by the intervention strategy \emph{regardless of the modelling uncertainty}. This latter aspect could be very valuable to epidemiologists because the parameters appearing in epidemic models sometimes cannot be accurately estimated (see \eg~\cite{Lobry2021}).

The contributions of the paper are as follows:
\begin{itemize}
	\item To our knowledge, this is the first paper that deals with the problem of maintaining a hard infection cap for the SIR and SEIR epidemic models (under the assumption of both perfect and imperfect modelling), via a set-theoretic approach: the theory of barriers, \cite{DonL13,EstAS20}.
	\item With the theory of barriers we are not only able to decrease the dimension in the computation of the sets, but also able to obtain easily checkable inequalities of the systems' parameters that allow one to classify the two sets, for the considered models, into trivial or nontrivial ones.	
	\item We demonstrate that the MRPI, whose use in epidemics was introduced in our previous research \cite{EstALetal20}, is also useful as a tool to maintain an infection cap for uncertain epidemic models.
	\item A note-worthy observation made in the paper regards the basic reproduction number, $\mathcal{R}_0$. In epidemiology intervention strategies are often designed so that $\mathcal{R}_0$ is less than one, which guarantees that the disease will eventually die out. We demonstrate that using the computed sets the infection cap can be maintained regardless of the value of $\mathcal{R}_0$, even if it is greater than one.
\end{itemize}

The paper is organised as follows. We present the SIR and SEIR models we will study in Section~\ref{sec_two_models}, along with precise formulations of the three problems we will address in the paper. We summarise the relevant results from the theory of barriers in Section~\ref{sec:theory_of_barriers}, which we apply to describe the sets for the two epidemic models in Sections~\ref{SIR_analysis}, \ref{section_MRPI_perfect_models} and \ref{section_MRPI_imperfect_models}. In Section~\ref{subsec_curves_backwards} we present conditions involving the system parameters that allow one to easily determine whether the sets are trivial or not. Detailed examples are presented in Section~\ref{Sec_examples} and the paper is concluded with a discussion in Section~\ref{sec_perspectives}.

\section{Two Epidemic Models}\label{sec_two_models}

In this section we describe the two well-known epidemic models on which the paper will focus: the SIR and SEIR models.

\subsection{SIR Model}\label{sec:SIR}

The (normalised) SIR model, see for example \cite[Ch. 2]{Brauer2019} and \cite{Het2000}, is described by the following ordinary differential equations:
\begin{align}
\dot S(t) &  =  - \beta S(t) I(t),\label{eq_SIR_nominal_1}\\
\dot I(t) &  =  \beta S(t) I(t) - \gamma I(t), \label{eq_SIR_nominal_2}\\
\dot R(t) & = \gamma I(t),\label{eq_SIR_nominal_3}
\end{align}
with $t\in[t_0,\infty[$ denoting time, measured in \emph{days}. The state variable $S(t)\in[0,1]$ denotes the proportion of individuals at time $t$ that are \emph{susceptible} to the disease; $I(t)\in[0,1]$ denotes the proportion that are infected with the disease and are \emph{infective} (that is, can infect susceptibles); and $R(t)\in[0,1]$ denotes the proportion that have been infected and \emph{removed}, in the sense that they cannot re-infect susceptibles. ``Removal'' may encapsulate recovery and immunisation from the disease, death from the disease, quarantine measures, vaccination, etc.

We assume that the birth and death rates are zero. In other words, we assume that the disease is relatively short-lived such that the effects of demographic change can be ignored. Therefore, assuming an initial state satisfying $S(t_0) + I(t_0) + R(t_0) = 1$, all integral curves of the system remain in the set: 
\begin{equation}
	\Pi = \{(S,I,R)^\top\in [0,1]^3 : S + I + R = 1\}, \label{Pi_SIR}
\end{equation}
for all $t\geq t_0$. Let $N\in\mathbb{R}_{>0}$ denote the constant population size, and let $s\triangleq NS$, $i\triangleq NI$ and $r\triangleq NR$ denote the number of susceptibles, infectives and recovered individuals, respectively. Let $\beta\in\mathbb{R}_{> 0}$, the \emph{contact rate}, denote the average number of contacts a person makes per unit time that is sufficient for disease transmission (this contact can be with infectives and/or other susceptibles). Then $\beta I$ is the average number of contacts with infectives a susceptible makes per unit time, and $\beta I s$ is the number of new infectives per unit time. Thus, $\beta I S$ is the \emph{proportion} of new infectives per unit time. It is assumed that recoveries occur according to a Poisson process with arrival rate $\gamma\in\mathbb{R}_{>0}$. See \cite[Ch. 2]{Brauer2019}, \cite{Het2000} and \cite{grundel2020much} for a more in-depth discussion of this model's derivation.

In conclusion, $\beta$ may be interpreted as the average number of people a member of the population comes into contact with over one day, and $\frac{1}{\gamma}$ is the average number of days an infective individual remains infective.

\subsection{SEIR Model}\label{sec:SEIR}

This model is described by the following ODEs:
\begin{align}
\dot S(t) &  =  - \beta S(t) I(t), \label{eq_SEIR_model_1} \\
\dot E(t) &  =  \beta S(t) I(t) - \eta E(t),\label{eq_SEIR_model_2}\\
\dot I(t) &  =  \eta E(t) - \gamma I(t),\label{eq_SEIR_model_3}\\
\dot R(t) & =   \gamma I(t),\label{eq_SEIR_model_4}
\end{align}
$t\in [t_0,\infty[$, with the additional compartment $E(t)\in[0,1]$ denoting \emph{exposed} individuals. This is a general model for diseases where individuals only become infective after a non-negligible period of time from exposure. It is assumed that an individual's evolution from being exposed to being infective follows a Poisson process with arrival rate $\eta\in\mathbb{R}_{>0}$. Thus, $\frac{1}{\eta}$ is the average number of days it takes an exposed individual to become infective.

As before, we assume that the total population remains constant and that the SEIR model is normalised. Thus, all solutions remain in the (redefined) set 
\begin{equation}
	\Pi =\{(S,E,I,R)^{\top}\in[0,1]^4 : S+E+I+R = 1\} \label{Pi_SEIR}
\end{equation}
for all $t\in [t_0,\infty[$.

\subsection{On Set-Based Methods Applied to Epidemiology}\label{Sec_use_of_sets}

Now we present an introductory overview of how we intend to apply set-based methods to the management of epidemics. The subsequent sections will present rigorous treatment of the ideas.

Consider the SIR model, \eqref{eq_SIR_nominal_1}-\eqref{eq_SIR_nominal_3}, and impose a hard infection cap, $I(t)\leq I_{\max}$, for all time. Moreover, assume that the contact rate, $\beta$, may be changed over time within some bounds, representing social distancing measures. That is, assume $\beta(t)\in[\beta_{\min},\beta_{\max}]$, $0 < \beta_{\min}\leq \beta_{\max}$ where $\beta_{\min}$ is the \emph{minimal} contact rate (and thus represents the \emph{maximal} social distancing) and $\beta_{\max}$ is the \emph{maximal} or \emph{nominal} contact rate (and may thus represent \emph{minimal} social distancing).

Consider the SEIR model, \eqref{eq_SEIR_model_1}-\eqref{eq_SEIR_model_4}, and also impose a hard infection cap, $I(t)\leq I_{\max}$, for all time. Assume that in addition to the contact rate, $\beta$, the parameter $\gamma$ may also be changed over time within some bounds, i.e., $\gamma(t)\in[\gamma_{\min}, \gamma_{\max}], 0<\gamma_{\min}\leq\gamma_{\max}$. This may represent quarantine or vaccination measures.

The first problem we want to address is the following:
\begin{problem}\label{Pb1:def}
For the SIR model, supposing the parameter $\gamma$ is perfectly known (SIR prefect model), how should the contact rate $\beta(t)$, be specified over time such that the infection cap, $I_{\max}$, is never breached? 
For the SEIR model, supposing the parameter $\eta$ is perfectly known (SEIR prefect model), how should the contact rate $\beta(t)$ and the removal rate $\gamma(t)$ be specified over time such that the infection cap is never breached?
\end{problem}
Our approach to solving this problem will be to find the \emph{admissible set}, $\A$, which, for the SIR model, will correspond to all states from which \emph{there exists an input} $\beta$ satisfying $\beta(t)\in[\beta_{\min},\beta_{\max}]$ for all  $t\in [t_0,\infty[$ and such that the infection cap, $I_{\max}$, is never breached. For the SEIR model, $\A$ will correspond to all states for which inputs $\beta(t)\in[\beta_{\min},\beta_{\max}]$ and $\gamma(t) \in[\gamma_{\min},\gamma_{\max}]$ exist such that the infection cap is never breached. Then the contact rate $\beta$, as well as the removal rate $\gamma$ for the SEIR model, will be deduced depending on the location of the state with respect to the computed sets. Moreover, the complement of the admissible set, $\Acomp$, is the set of states for which the infection cap \emph{cannot be maintained}, regardless of the interventions. Thus, from this set the problem has no solution and if the state is located here, either the infection cap needs to be relaxed, or the maximal allowable social distancing measures ($\beta_{\min}$), and/or maximal quarantine rate for the SEIR model ($\gamma_{\max}$), need to be strengthened. 

The set $\A$ is obtained by constructing its \emph{boundary}, which will be seen to be made of two parts, the \emph{usable part}, a subset included in the boundary of the state constraints, and the \emph{barrier}, entirely contained in the interior of the state constraint set (see the general presentation of Section~\ref{SIR_analysis},  and further results and comparisons for the SIR and SEIR perfect models in Sections~\ref{subsec_curves_backwards} and \ref{Sec_examples}).


The second problem concerns also the SIR and SEIR perfect models:
\begin{problem}\label{Pb2:def}
For the SIR perfect model, we want to determine the set of initial states, $(S(t_0),I(t_0),R(t_0))$, from which the infection cap $I_{\max}$ will be maintained for all times and for any function $\beta: [t_0,\infty[\rightarrow [\beta_{\min}, \beta_{\max}]$, $\gamma$ being perfectly known. 
For the SEIR model, the problem consists also in maintaining the infection cap $I_{\max}$ for all times and any pair of functions $(\beta,\gamma)$ with $\beta: [t_0,\infty[\rightarrow [\beta_{\min}, \beta_{\max}]$ and  $\gamma: [t_0,\infty[\rightarrow [\gamma_{\min}, \gamma_{\max}]$, $\eta$ being perfectly known.  \end{problem}
 These sets are called the \emph{Maximal Robust Positively Invariant (MRPI)} sets associated to the SIR and SEIR perfect model respectively (see their general definitions~\ref{RPI:def} and \ref{MRPI:def}). Their construction, again via the construction of their boundary, consisting, as before, of a usable part and a barrier, is presented in Sections~\ref{section_MRPI_perfect_models}, \ref{subsec_curves_backwards} and \ref{Sec_examples}).

The third problem concerns the SIR and SEIR imperfect models:
\begin{problem}\label{Pb3:def}
Suppose now that, for the SIR model,  the value of the parameter $\gamma$ is \emph{not} known. All that we know is that $\gamma\in [\gamma_{\min}, \gamma_{\max}]$. We also assume that an intervention strategy $\hat{\beta}(t):[t_0,\infty[\rightarrow [\beta_{\min},\beta_{\max}]$ has been \emph{pre-designed} (SIR imperfect model). Find the states from which this intervention strategy  will maintain the infection cap \emph{regardless of the error in the parameter $\gamma$}, for all $t\in [t_0,\infty[$. 
For the SEIR model, we pose the same problem with unknown parameter $\eta \in [\eta_{\min}, \eta_{\max}]$  in place of $\gamma$ and with predesigned intervention strategies $\hat{\beta}(t):[t_0,\infty[\rightarrow [\beta_{\min},\beta_{\max}]$ and $\hat{\gamma}(t):[t_0,\infty[\rightarrow [\gamma_{\min},\gamma_{\max}]$ (SEIR imperfect model).
 \end{problem}
 
Its solution will consist in constructing the MRPI set associated to the SIR (resp. SEIR) imperfect model (see Sections~\ref{section_MRPI_imperfect_models}, \ref{subsec_curves_backwards} and \ref{Sec_examples}). 

Table~\ref{table} summarises the status of the parameters and controls for each problem.
\begin{center}
\begin{table}[h]\label{table}
\caption{Summary of known/unknown parameters and controlled/uncontrolled inputs.}\label{table}
 \begin{tabular}{||c|c||c|c|c|c||}
 \hline
 Problem&Model&known&unknown&controllable&pre-designed\\
 &&param.&param.&input&feedback\\
 \hline
1 and 2&SIR Perfect& $\gamma$ &- & $\beta$ & -\\
 \hline
1 and 2&SEIR Perfect& $\eta$ & - & $\beta, \gamma$ & -\\
 \hline
3&SIR Imperfect & - & $\gamma$ & - & $\hat{\beta}$\\
 \hline
3&SEIR Imperfect & - & $\eta$ & - & $\hat{\beta}, \hat{\gamma}$\\
 \hline
\end{tabular}
\end{table}
\end{center}
The main mathematical tool to solve these three problems is sketched in the next section.

\section{The Theory of Barriers}\label{sec:theory_of_barriers}

We now summarise the relevant theory from \cite{DonL13} and \cite{EstAS20}, which we will apply to describe the sets of the two epidemic models. Consider a state and input constrained nonlinear system:
\begin{align}
\dot{x}(t) & = f(x(t),u(t)),\quad x(t_0) = x_0,\quad u\in\mathcal{U},\label{sys_eq_1}\\
g_i(x(t)) &\leq 0,\forall t\in[t_0,\infty[,\,\,i=1,2,\dots,p,\label{sys_eq_2}
\end{align}
where $x(t)\in\mathbb{R}^n$ is the state; $x_0$ is the initial condition; $u(t)\in\mathbb{R}^m$ is the input; $\mathcal{U}$ is the set of Lebesgue measurable functions that map the interval $[t_0,\infty[\subset\mathbb{R}$ to a compact and convex set $U\subset\mathbb{R}^m$ (note the difference between $\mathcal{U}$, a function space, and $U$, a subset of $\RR^m$); and $g_i$, $i=1,\dots,p$, are the state constraints. We impose the following assumptions:
\begin{description}
	\item[Assumptions]
\end{description}
\begin{itemize}
	\item[(A1)] The function $f:\mathbb{R}^n\times\mathbb{R}^m\rightarrow\mathbb{R}^n$ is $\mathcal{C}^2$ on an open set containing $\mathbb{R}^n \times U$.
	\item[(A2)] Every $x_0$, with $\Vert x_0 \Vert < \infty$, and every $u\in\mathcal{U}$ admits a unique absolutely continuous integral curve of \eqref{sys_eq_1} that remains bounded over any finite interval of time.
	\item[(A3)] The set $f(x,U)$ is convex for all $x\in\mathbb{R}^n$.
	\item[(A4)] The function $g_i$, $i=1,\dots,p$, is $\mathcal{C}^2$ and the set of points given by $g_i(x) = 0$ defines an $n-1$ dimensional manifold.
\end{itemize}
It can be verified that the SIR and SEIR models analysed in subsequent sections satisfy Assumptions (A1)-(A4) for each considered interpretation of the input. In particular, the following growth condition is always satisfied: there exists a $C<\infty$ such that $\sup_{u\in \U} |x^\top f(x,u)| \leq C (1 + \Vert x \Vert^2)$ for all $x\in\mathbb{R}^n$, which implies (A2). An important consequence of these assumptions is that the admissible and MRPI sets are closed.

We use the following notation: let $x^{(u,x_0,t_0)}$ denote the unique solution of \eqref{sys_eq_1} from the initial state $x_0\in\RR^n$ at initial time $t_0$, with input function $u\in \U$, and let $x^{(u,x_0,t_0)}(t)$ denote the solution at time $t\in[t_0,\infty[$. If clear from context, we will use the notation $x^{u}$ and $x^{u}(t)$ instead. Let $G\triangleq \{x : g_i(x)\leq 0, i = 1,2,...,p\}$, $G_-\triangleq \{x : g_i(x) <0, i = 1,2,...,p\}$, and $G_0\triangleq \{x : \exists i\in\{1,2,...,p\} \;\; \text{s.t.} \; g_i(x) =0\}.$ The set of indices of active constraints at a point $x$ is denoted $\mathbb{I}(x)=\{i: g_i(x) = 0\}$. We let $L_fg(x,u)\triangleq Dg(x)f(x,u)$ denote the Lie derivative of a differentiable function $g:\mathbb{R}^n \rightarrow \mathbb{R}$ with respect to $f(\cdot,u)$ at the point $x$.

We now present the definitions of the two sets that we will describe for the compartmental epidemic models.
\begin{definition}
	The \emph{admissible set}\footnote{also known as the \emph{viability kernel} in Viability theory, \cite{aubin2009viability}}, see \cite{DonL13}, of the system \eqref{sys_eq_1}-\eqref{sys_eq_2}, denoted by $\A$, is the set of initial states for which there exists a $u\in \U$ such that the corresponding solution satisfies the constraints \eqref{sys_eq_2} for all future time,
	\begin{align*}
	\mathcal{A} \triangleq \left\{ x_0 \in G: \exists u\in\mathcal{U}, \;\;  x^{(u,x_0,t_0)}(t)\in G \;\; \forall t \in [t_0,\infty[ \;\right\}.
	\end{align*}
\end{definition}
\begin{definition}\label{RPI:def}
	A set $\Omega\subset\RR^n$ is said to be a \emph{robust positively invariant set} (RPI) of the system \eqref{sys_eq_1} provided that $x^{(u,x_0,t_0)}(t)\in\Omega$ for all $t\in[t_0,\infty[$, for all $x_0\in\Omega$ and for all $u\in\U$.
\end{definition}
\begin{definition}\label{MRPI:def}
	The \emph{maximal robust positively invariant set} (MRPI) of the system \eqref{sys_eq_1}-\eqref{sys_eq_2} contained in $G$, see \cite{EstAS20}, is the union of all RPIs that are subsets of $G$. Equivalently\footnote{this equivalence is shown in \cite{EstAS20}},
	\begin{align*}
	\mathcal{M} \triangleq \left\{ x_0\in G: x^{(u,x_0,t_0)}(t)\in G, \,\,\forall u\in\mathcal{U}, \;\;  \forall t \in [t_0,\infty[ \;\right\}.
	\end{align*}
\end{definition}
We obviously have $\mathcal{M}  \subset \A$.
Moreover, under (A1)-(A4) the sets $\A$ and $\M$ are closed, and we denote their boundaries by $\partial \A$ and $\partial \M$, respectively.

The main result from the papers \cite{DonL13} and \cite{EstAS20} is a characterisation of the sets' boundaries that are made up of two complementary parts: the \emph{usable part}, which is contained in the boundary of the constrained state space, $\DAO \triangleq  \partial\A\cap G_0$ and $\DMO \triangleq  \partial\M\cap G_0$; and the \emph{barrier}, which is contained in the interior of the constrained state space, $\DAM \triangleq  \partial\A\cap G_-$ and $\DMM \triangleq  \partial\M\cap G_-$. Points on the admissible set's usable part satisfy:
\begin{equation}\label{A-usable_part_eq}
\min_{u\in U} \max_{i\in\mathbb{I}(x)} L_f g_i (x,u) \leq 0, \,\,\forall x\in\DAO,
\end{equation}
whereas those on the MRPI's usable part satisfy:
\begin{equation}\label{MRPI-usable_part_eq}
\max_{u\in U} \max_{i\in\mathbb{I}(x)} L_f g_i (x,u) \leq 0, \,\,\forall x\in\DMO.
\end{equation}
The barrier is made up of integral curves of the system that satisfy a minimum-/maximum-like principle that may intersect the boundary of the constraint set, $G_0$, in a tangential manner. We summarise these results in the following theorem, for which proofs can be found in \cite{DonL13} and \cite{EstAS20}.
\begin{theorem}\label{thm1}
	Under the assumptions (A1)-(A4),  every integral curve $x^{\bar{u}}$ contained in $\DAM$ (resp. $\DMM$) and crossing $G_0$ tangentially in finite time satisfy, together with the corresponding control $\bar{u}\in\U$, the following necessary conditions. There exists a nonzero absolutely continuous maximal solution $\lambda^{\bar{u}}$ to the adjoint equation:
	\begin{equation}
	\dot{\lambda}^{\bar{u}}(t) = -\left( \frac{\partial f}{\partial x}(x^{\bar{u}}(t),\bar{u}(t)) \right)^T \lambda^{\bar{u}}(t),\label{eq_adjoint}
	\end{equation}
	such that the \emph{Hamiltonian} $\lambda^T f$ satisfies:
	\begin{align}		
	\min_{u\in U}\{\lambda^{\bar{u}}(t)^Tf(x^{\bar{u}}(t),u) \}&= \lambda^{\bar{u}}(t)^T f(x^{\bar{u}}(t),\bar{u}(t)) = 0, \label{eq_Hamil_A} \\
	\Big(\Big.\text{resp.}\;
	\max_{u\in U}\{\lambda^{\bar{u}}(t)^Tf(x^{\bar{u}}(t),u) \} &= \lambda^{\bar{u}}(t)^T f(x^{\bar{u}}(t),\bar{u}(t)) = 0\label{eq_Ham_max_M}
	\Big.\Big),
	\end{align}
	for almost all $t$. Moreover, we have:
	$\lambda^{\bar{u}}(\bar{t}) = (Dg_{i^{\star}}(z))^T,$ 
	where
	\begin{align}
	&\min_{u\in U}\max_{i\in\mathbb{I}(z)} L_fg_i(z,u) = L_fg_{i^{\star}}(z,\bar{u}(\bar{t})) = 0, \label{thm1_ult_tan_A} \\
	\Big(\Big.\text{resp.}\quad
	&\max_{u\in U}\max_{i\in\mathbb{I}(z)} L_fg_i(z,u) = L_fg_{i^{\star}}(z,\bar{u}(\bar{t})) = 0 \label{thm1_ult_tan_M}
	\Big.\Big),
	\end{align}
	$\bar{t}$ denotes the time at which $G_0$ is reached, $z \triangleq x^{\bar{u}}(\bar{t})\in G_0$, and $i^{\star}\in \{1, \ldots, p\}$ denotes any index for which the constraints are active at the point $z$, i.e. $i^{\star} \in \mathbb{I}(z)$.
\end{theorem}
To now compute the curves running along the barrier, we first obtain the points of \emph{ultimate tangentiality} via the expressions \eqref{thm1_ult_tan_A} (resp. \eqref{thm1_ult_tan_M}) and then integrate the system backwards in time with the control function that minimizes (in the case of the set $\A$) or maximizes (in the case of the MRPI set $\M$) the \emph{Hamiltonian} for almost every $t$, according to the expression \eqref{eq_Hamil_A} (resp. \eqref{eq_Ham_max_M}) with $\lambda$ satisfying \eqref{eq_adjoint}.

In the sequel the set $G$ will be expressed by a single state constraint that caps the proportion of infective individuals, and we will interpret the function $u$ in two ways. The first interpretation will be that $u$ is a time-varying ``controllable'' input ($u=\beta$ for the SIR model and $u=(\beta,\gamma)$ for the SEIR model), which will allow us to solve Problem~\ref{Pb1:def} (resp. Problem~\ref{Pb2:def}) (see Subsection~\ref{Sec_use_of_sets} and Table~\ref{table}), \ie, find the contact rate $\beta$ (in the SEIR case also the removal rate $\gamma$) such that, the parameter $\gamma$ being perfectly known, the infection cap, $I_{\max}$, is never breached (resp. never breached for all $\beta \in [\beta_{\min},\beta_{\max}]$).
	
The second interpretation will be that $u$ is an ``uncontrollable'' \emph{disturbance} term ($u=\gamma$ for the SIR model and $u=\eta$ for the SEIR model) encapsulating uncertainty in modelling parameters. This will allow us to address Problem~\ref{Pb3:def}, \ie, given a pre-designed feedback, find the states from which the infection cap will be always maintained  \emph{regardless of the error in the modelling parameters}.

\section{Carrying out the Analysis of the Admissible Set for the Perfect SIR and SEIR Models (Problem~\ref{Pb1:def})}\label{SIR_analysis}

In this section we carry out the detailed analysis of the admissible sets of the constrained SIR and SEIR models under the assumption of perfect modelling (Problem~\ref{Pb1:def}).
For the sake of readability, comparisons of admissible and robust sets for the perfect and imperfect models as well as numerical applications are presented in Sections~\ref{subsec_curves_backwards} and \ref{Sec_examples}.

\subsection{Admissible Set for Perfect SIR Model}\label{subsection:perfect_model_case}

We assume that the SIR model is perfect (that is, $\gamma$ is known perfectly) and that $\beta$ is a controllable time-varying input, which may, for example, model social-distancing measures to halt the disease's spread. We impose hard constraints on the input and state: $I(t)\leq I_{\max}$ and $\beta(t)\in[\beta_{\min},\beta_{\max}]$ for all $t\in[t_0,\infty[$, where $I_{\max}\in[0,1]$ is the hard infection cap, and $\beta_{\min},\beta_{\max}\in\mathbb{R}_{>0}$, with $0<\beta_{\min}\leq \beta_{\max}$, are the minimal and maximal contact rates, respectively. Recall that the bound $\beta_{\max}$ may be interpreted as the \emph{nominal contact rate} in absence of the disease, and $\beta_{\min}$ as the \emph{maximal social distancing}.

Because $R$ has no influence on $S$ or $I$ we can study the two-dimensional subsystem:
\begin{align}
\dot S(t) &  =  - \beta(t) S(t) I(t), \label{eq_SIR_1}\\
\dot I(t) &  =  \beta(t) S(t) I(t) - \gamma I(t), \label{eq_SIR_2}\\
I(t) - I_{\max} & \leq 0,\quad \beta(t) \in [\beta_{\min},\beta_{\max}],\label{eq_SIR_3}
\end{align}
with $t\in[t_0,\infty[$. 
To ease our notation we will sometimes label the state $x\triangleq (S,I)^\top \in \mathbb{R}^2$, the control $u\triangleq\beta\in\mathbb{R}$, the control constraint set $U \triangleq [\beta_{\min}, \beta_{\max}]\subset \mathbb{R}$, the single constraint function $g(x) \triangleq I - I_{\max}$, and the set $G_{\Pi}\triangleq G\cap\Pi$, where $\Pi$ is given by \eqref{Pi_SIR}. We also denote by $f(x,\beta)= 
\left( \begin{array}{c} -\beta SI\\
\beta SI - \gamma I
\end{array}\right)$.

Applying condition \eqref{A-usable_part_eq}, the usable part $\DAO$ is given by the condition
$$\min_{\beta \in [\beta_{\min}, \beta_{\max}]} L_f g = \beta SI - \gamma I \leq 0, \quad (S,I)\in G_{\Pi}~\mathrm{s.t.~} I-I_{\max}=0,$$
which yields $(\beta_{\min} S- \gamma) I_{\max} \leq 0$, or $S\leq \frac{\gamma}{\beta_{\min}}$.
Therefore the usable part is the set $\DAO = \{(S,I)\in G_{\Pi} \mid 0\leq S \leq \min \left( \frac{\gamma}{\beta_{\min}}, 1-I_{\max}\right), \; I=I_{\max}\}$.
\bigskip

Focusing now on the barrier $\DAM$, let us label tangent points to $G_{0}$ by $z\triangleq (z_1,z_2)^\top = (S(\bar{t}),I(\bar{t}))^\top \in \mathbb{R}^2$, where $\bar{t}$ denotes the time at which $G_{0}$ is reached, as explained in Theorem~\ref{thm1}. Invoking ultimate tangentiality, \eqref{thm1_ult_tan_A}, we get: 
$Dg(z_1,z_2)= (0,1)$ and $f(z,\beta)= 
\left( \begin{array}{c} -\beta z_{1}z_{2}\\
\beta z_{1}z_{2} - \gamma z_{2}
\end{array}\right)$.
Thus:
\begin{equation}\label{tancondSIRp:eq}
\begin{aligned}
		\min_{u\in U} L_f g(z,u) &=  \min_{\beta \in [\beta_{\min}, \beta_{\max}]} Dg(z_{1},z_{2})f(z_{1},z_{2},\beta) \\
		&= \min_{\beta \in [\beta_{\min}, \beta_{\max}]} \{\beta z_2z_1 - \gamma z_2\} = 0,
\end{aligned}
\end{equation}
and along with the constraints $z_2 = I_{\max}$ and $z_1 + z_2 \leq 1$, we get that 
the single tangent point must satisfy 
$(\beta_{\min}z_{1} -\gamma)  I_{\max} =0$, or $(z_1,z_2) = (\frac{\gamma}{\beta_{\min}}, I_{\max})$, with $ \frac{\gamma}{\beta_{\min}} + I_{\max} \leq 1$, one of the end point of the usable part. 

The adjoint system, \eqref{eq_adjoint}, reads:
\begin{equation}
\dot{\lambda}(t) = 
\left(
\begin{array}{ccc}
\beta(t) I(t) & - \beta(t) I(t) \\
\beta(t) S(t) & - \beta(t) S(t)  + \gamma
\end{array}
\right)\lambda(t),\,\,\lambda(\bar{t}) = (0,1)^T\label{SIR_perfect_adjoint}
\end{equation}
From the Hamiltonian minimisation condition, \eqref{eq_Hamil_A}, we get:
\begin{equation}\label{SIR_Hamil_min:eq}
\min_{\beta(t)\in [\beta_{\min},\beta_{\max}]} \lambda_1(t)\left[ -\beta(t)S(t)I(t) \right] + \lambda_2(t)\left[\beta(t)S(t)I(t) - \gamma I(t) \right] = 0,
\end{equation}
from where we get:
\begin{equation}
\bar{\beta}(t) = 
\begin{cases} 
\beta_{\max}&\mbox{if } \lambda_2(t) - \lambda_1(t) < 0, \\
\beta_{\min}&\mbox{if } \lambda_2(t) - \lambda_1(t) > 0, \\
\textrm{anything}&\mbox{if } \lambda_2(t) - \lambda_1(t) = 0,
\end{cases}\label{input_SIR_perfect_admissible}
\end{equation}
where $\bar{\beta}$ denotes the input associated with the barrier trajectories. Thus, at the point $z$, we have $\lambda_{2}(\bar{t})-\lambda_{1}(\bar{t})=1$, which immediately implies, by continuity of the adjoint, that $\bar{\beta}(t)= \beta_{\min}$ in some nonempty interval, that is, for $t\in]\bar{t}-\ee, \bar{t}]$, $\ee >0$. 

\begin{proposition}\label{prop_SIR_perfect_beta_non_singular}
	The input $\bar{\beta}$, given by \eqref{input_SIR_perfect_admissible}, associated with the barrier of the admissible set of the perfect SIR model, \eqref{eq_SIR_1}-\eqref{eq_SIR_3} only switches at isolated points in time.	
\end{proposition}
\begin{proof}
	Suppose that $\lambda_2(t) - \lambda_1(t)= 0$ over an interval of time $]\hat{t} - \ee, \hat{t} + \ee[$, with $\hat{t}\in\mathbb{R}$ arbitrary and $\ee>0$. Then, we would have $\dot{\lambda}_2(t) - \dot{\lambda}_1(t) =0$ identically in this interval, \ie, according to \eqref{SIR_perfect_adjoint}, 
$$
\dot{\lambda}_2(t) - \dot{\lambda}_1(t) =	\beta(t)(\lambda_2(t) - \lambda_1(t))(I(t) - S(t)) + \gamma\lambda_2(t) = \gamma\lambda_2(t) = 0, \quad \forall t\in]\hat{t} - \ee, \hat{t} + \ee[.$$ 
We deduce that $\gamma \lambda_2(t)= \gamma \lambda_1(t) = 0$ over this interval. Therefore, the adjoint would identically vanish, which is impossible according to Theorem~\ref{thm1}, hence the proposition.
\end{proof}
\begin{remark}
	As is known from optimal control, an input with bang-bang structure may exhibit other complicated behaviour, such as the ``Fuller phenomenon''. Though we cannot make definite statements concerning the existence or nonexistence of such behaviours along barrier trajectories, they were not observed in any of the numerical examples in this paper.
\end{remark}

The next Proposition states that barrier curves associated with the admissible set of the perfect SIR model always evolve backwards into $G_-$.

\begin{proposition}\label{prop:backwards_evolve}
Consider the constrained perfect SIR model, \eqref{eq_SIR_1}-\eqref{eq_SIR_3}. The barrier $\DAM$ is constituted by the curve $x^{\bar{u}}$ generated by \eqref{input_SIR_perfect_admissible} with \eqref{SIR_perfect_adjoint}, that ends tangentially to $G_0\cap \Pi$ at the point $(\frac{\gamma}{\beta_{\min}}, I_{\max})\in G_0\cap \Pi$ with $\frac{\gamma}{\beta_{\min}} + I_{\max} \leq 1$, and is contained in $G_-$ at least for all $t\in]\bar{t} - \ee,\bar{t}[$, $\ee >0$.
\end{proposition}
\begin{proof}
	Let $\bar{g}$ denote the mapping $t\rightarrow \bar{g}(t)\triangleq g(x^{\bar{u}}(t)) = I(t) - I_{\max}$ over some open interval $]\bar{t} - \ee,\bar{t}[$, with $\ee >0$.  We have shown that the barrier curve $x^{\bar{u}}= (S(t),I(t))$ has to be such that $\bar{g}$ vanishes at $\bar{t}$ as well as its first derivative, which is equal to $\frac{d\bar{g}}{dt} = L_{f}g$. Moreover, $\beta$ must remain constant, equal to $\beta_{\min}$, in the interval $]\bar{t} - \ee,\bar{t}[$. Let us compute $\bar{g}$'s left second derivative:
	\begin{equation*}
	\begin{aligned}
		&\partial_- L_fg(x^{\bar{u}}(\bar{t}),\bar{u}(\bar{t})) \triangleq \lim_{t\rightarrow \bar{t}, \,t\leq \bar{t}} \frac{d}{dt} \left (\beta(t) I(t)S(t) - \gamma I(t)\right) \\
		&= \beta(\bar{t})I(\bar{t})\left[-\beta(\bar{t})S(\bar{t})I(\bar{t})\right] + \beta(\bar{t})\left[\beta(\bar{t})S(\bar{t})I(\bar{t})-\gamma I(\bar{t})\right]S(\bar{t}) - \gamma\left[\beta(\bar{t})S(\bar{t})I(\bar{t})-\gamma I(\bar{t})\right].
		\end{aligned}
	\end{equation*}
	Using the fact that $\beta(\bar{t})S(\bar{t})I(\bar{t})- \gamma I(\bar{t}) = 0$, $S(\bar{t}) = \frac{\gamma}{\beta_{\min}}$, $\beta(\bar{t})=\beta_{\min}$ and $I(\bar{t})=I_{\max}$, this simplifies to:
	\begin{equation*}
		\partial_- L_fg(x^{\bar{u}}(\bar{t}),\bar{u}(\bar{t})) = - \gamma \beta_{\min}I_{\max}^2 < 0.
	\end{equation*}
	Therefore, using Taylor's expansion formula, 
	$$\begin{aligned}
	g(x^{\bar{u}}(t)) &= g(x^{\bar{u}}(\bar{t})) + (t-\bar{t}) L_{f}g(x^{\bar{u}}(\bar{t})) + \frac{1}{2} (t-\bar{t})^2 \partial_- L_fg(x^{\bar{u}}(\bar{t}),\bar{u}(\bar{t})) + 0((t-\bar{t})^3) \\
	&= 0 + 0 - \frac{1}{2} (t-\bar{t})^2 \beta_{\min}I_{\max}^2\gamma + 0((t-\bar{t})^3) \leq 0
	\end{aligned}$$ 
	in an interval $]\bar{t} - \ee,\bar{t}[$ for a sufficiently small $\ee$, which proves that the barrier curve $x^{\bar{u}}$ remains in $G_-$ in this interval.
\end{proof}

\begin{remark}\label{barrier_extension:rem}
Note that the point of the barrier satisfying $I=0$ is an equilibrium point: $\dot{S}= -\beta SI=0$, $\dot{I}= (\beta S - \gamma)I=0$. Therefore the unique solution from $(S(t_{0}),0)$ is reduced to this point. 

In fact, the whole axis $I=0$ may be seen to be included in $\DAM$ since every point of this axis is an equilibrium and thus satisfies the constraint forever and belongs to $\partial \Pi$, the boundary of $\Pi$ (see also Section~\ref{subsecSIR_curves_backwards}). It is remarkable that the points in this axis, though in the barrier, do not satisfy the conditions of Theorem~\ref{thm1} since they do not cross $G_{0}$ (recall that Theorem~\ref{thm1} is only necessary).

Finally, the integral curve of the system $S(t)=0$, $I(t) = e^{-\gamma (\bar{t}-t)}I_{\max}$, $t\in ]-\infty, \bar{t}]$, \ie~ the semi-open line segment $\{(0, I) \mid I\in ]0,I_{\max}]\}$ also lies in $\A\cap \partial \Pi$ and hence in $\DAM$, neither satisfies the conditions of Theorem~\ref{thm1} since it does not intersect $G_{0}$ tangentially. It results that $\DAM$ is the union of the set made by the integral curves satisfying Theorem~\ref{thm1} and of a set of special curves that do not intersect $G_{0}$ tangentially.
\end{remark}
%

\subsection{Admissible Set for Perfect SEIR Model}\label{Subsection:admissible_set_SEIR}

We now focus on the model \eqref{eq_SEIR_model_1}-\eqref{eq_SEIR_model_4} and assume that $\beta$ and $\gamma$ are controllable inputs representing social distancing measures, and quarantining/isolation measures of infectives, respectively. We assume that $\eta$ is perfectly known. Thus, we consider the three-dimensional constrained system:
\begin{align}
\dot S(t) &  =  - \beta(t) S(t) I(t), \label{eq_SEIR_1}\\
\dot E(t) &  =  \beta(t) S(t) I(t) - \eta E(t), \label{eq_SEIR_2}\\
\dot I(t) &  =  \eta E(t) - \gamma(t) I(t), \label{eq_SEIR_3}\\
I(t) - I_{\max} & \leq 0,\quad \beta(t) \in [\beta_{\min},\beta_{\max}],\quad \gamma(t) \in [\gamma_{\min},\gamma_{\max}]\label{eq_SEIR_4}
\end{align}
with $t\in[t_0,\infty)$, $0< I_{\max}\leq 1$, $0<\beta_{\min}\leq \beta_{\max}$, and $0<\gamma_{\min}\leq \gamma_{\max}$. We will sometimes label the state $x\triangleq (S,E,I)^\top \in \mathbb{R}^3$, the controllable input $u\triangleq(\beta,\gamma)^{\top}\in\mathbb{R}^2$, the control constraint set $U \triangleq [\beta_{\min}, \beta_{\max}]\times[\gamma_{\min}, \gamma_{\max}]\subset \mathbb{R}^2$, $g(x) \triangleq I - I_{\max}$, and $G_{\Pi}\triangleq G \cap \Pi$, where $\Pi$ is given by \eqref{Pi_SEIR}.


As we did with the perfect SIR model, we first construct the usable part $\DAO$ given, using  \eqref{A-usable_part_eq}, by the condition
$$\min_{\gamma\in [\gamma_{\min}, \gamma_{\max}]} L_f g = \eta E - \gamma I \leq 0, \quad (S,E,I)\in G_{\Pi}~\mathrm{s.t.~} I-I_{\max}=0$$
which yields $\gamma =\gamma_{\max}$ and $E \leq \min\{\frac{\gamma_{\max}}{\eta}I_{\max} , \; 1 - S - I_{\max} \}$, $S\in [0,1-I_{\max}]$, or
$$\DAO = \{(S,E,I) \mid S  \in [0,1-I_{\max}], E \in [0, \min\{\frac{\gamma_{\max}}{\eta}I_{\max} , \; 1 - S - I_{\max}\}], I = I_{\max}\}.$$
\bigskip

Now, the barrier curves of $\DAM$ must satisfy the ultimate tangentiality condition, \eqref{thm1_ult_tan_A}, to get:
\begin{align*}
\min_{u\in U} L_f g(z,u) = \min_{\gamma\in [\gamma_{\min}, \gamma_{\max}]} \{\eta z_2 - \gamma z_3 \} = 0,
\end{align*}
and along with the constraints $z_3 = I_{\max}$ and $0\leq z_1 + z_2 + z_3 \leq 1$, we obtain the \emph{set} of tangent points:
\begin{equation}
\mathcal{T}_{\A}\triangleq\{ (z_1,z_2,z_3) \mid 0\leq z_1\leq 1 - z_2 - z_3; \, z_2 = \frac{\gamma_{\max}}{\eta} I_{\max}; \,z_3= I_{\max}\}.\label{eq_tan_pts_SEIR}
\end{equation}
Thus, unlike with the perfect SIR model, we now have a line segment in the boundary of $\DAO$ where the admissible set's barrier curves intersect $G_0$ tangentially. The adjoint system, \eqref{eq_adjoint}, reads:
\begin{equation}
\dot{\lambda} = 
\left(
\begin{array}{ccc}
\beta I & - \beta I & 0 \\
0 & \eta & -\eta \\
\beta S & - \beta S & \gamma
\end{array}
\right)\lambda,\,\,\lambda(\bar{t}) = (0,0,1)^T,\label{adjoint_SEIR_perfect}
\end{equation}
and condition~\eqref{eq_Hamil_A} gives the input associated with the admissible set's barrier curves:
\begin{equation}\label{mini-beta-gamma_SEIR:eq}
\bar{\beta}(t) = 
\begin{cases} 
\beta_{\max}&\mbox{if } \lambda_2(t) - \lambda_1(t) < 0, \\
\beta_{\min}&\mbox{if } \lambda_2(t) - \lambda_1(t) > 0,\\
\textrm{anything}&\mbox{if } \lambda_2(t) - \lambda_1(t) = 0,
\end{cases}
~\mathrm{and~}
\bar{\gamma}(t) = 
\begin{cases} 
\gamma_{\max}&\mbox{if } \lambda_3(t) >  0, \\
\gamma_{\min}&\mbox{if } \lambda_3(t) < 0,\\
\textrm{anything}&\mbox{if } \lambda_3(t) = 0.
\end{cases}
\end{equation}
Since $\lambda_{3}(\bar{t})=1$, we immediately see that $\bar{\gamma}(\bar{t}) = \gamma_{\max}$. Moreover, we see from \eqref{adjoint_SEIR_perfect} that $\dot{\lambda}_2(t) - \dot{\lambda}_1(t) = \eta(\lambda_{2}(t)-\lambda_{3}(t)) - \beta(t)I(t)(\lambda_{1}(t)-\lambda_{2}(t))$ and thus, at $\bar{t}$, we get
$\dot{\lambda}_2(\bar{t}) - \dot{\lambda}_1(\bar{t}) = -\eta < 0$, implying that $\bar{\beta}(t) = \beta_{\min}$ for all $t\in ]\bar{t} - \ee, \bar{t}]$, $\ee>0$. 

We have the following results concerning integral barrier curves.
\begin{proposition}\label{prop_SEIR_perfect_inputs_nonsingular}
Consider the constrained perfect SEIR model, \eqref{eq_SEIR_1}-\eqref{eq_SEIR_4}. Along any barrier integral curve satisfying \eqref{adjoint_SEIR_perfect}-\eqref{mini-beta-gamma_SEIR:eq} and such that $z_{1} \neq 0$ in \eqref{eq_tan_pts_SEIR}, the inputs $\bar{\beta}$ and $\bar{\gamma}$, given by \eqref{mini-beta-gamma_SEIR:eq}, only switch at isolated points in time.
\end{proposition}
\begin{proof}
	Following the same arguments as in Proposition~\ref{prop_SIR_perfect_beta_non_singular}, suppose $\lambda_2(t) - \lambda_1(t) = 0$ over some interval of time, $\mathcal{I}\triangleq ]\hat{t}-\ee, \hat{t} + \ee[$, $\hat{t}\in\mathbb{R}$, $\ee\in\mathbb{R}_{\geq 0}$. Then, we would have  $\dot{\lambda}_2(t) - \dot{\lambda}_1(t) = \eta(\lambda_2(t) - \lambda_3(t)) + \beta(t)I(t)(\lambda_2(t) - \lambda_1(t)) = 0$ for all $t\in \mathcal{I}$, which would imply $\lambda_2(t) - \lambda_3(t) = 0$ for all $t\in \mathcal{I}$.
This would also imply that $\lambda_3(t)=0$ for all $t\in \mathcal{I}$ since
$\dot{\lambda}_2(t) - \dot{\lambda_3}(t) = 0 = \eta(\lambda_2(t) - \lambda_3(t)) + \beta(t)S(t)(\lambda_2(t) - \lambda_1(t))   -\gamma(t)\lambda_{3}(t)    = -\gamma \lambda_3(t)$ for all $t\in \mathcal{I}$.
Moreover, noting that $\dot{\lambda}_3(t) = -\beta(t)S(t)(\lambda_2(t) - \lambda_1(t)) + \gamma\lambda_3(t)$, we would have $\dot{\lambda}_3(t) = \gamma\lambda_3(t) = 0$ for all $t\in I$.
Thus, $\lambda_1(t) = \lambda_2(t) = \lambda_3(t) = 0$ for all $t\in \mathcal{I}$, which is impossible according to Theorem~\ref{thm1}.
	
	Suppose now that $\lambda_3(t) = 0$ over some nonempty interval $\mathcal{I}$. Then we would have $\dot{\lambda}_3(t) = -\beta(t)S(t)(\lambda_2(t) - \lambda_1(t)) = 0$ over this interval, which would imply that $S(t)(\lambda_2(t) - \lambda_1(t)) = 0$, provided that $S(t)\neq 0$, over this interval. But if $S(t)= 0$ over some interval of time, according to the uniqueness of the solution of \eqref{eq_SEIR_1}-\eqref{eq_SEIR_4}, this would imply that $S(t)=0$ for all times and in particular $S(\bar{t})= z_{1}=0$, which is excluded by assumption.  Therefore $\lambda_3(t) = 0$ over some nonempty interval would imply $\lambda_2(t) - \lambda_1(t) = 0$ over this interval, which we have established is impossible in the first part of the proof.
\end{proof}

The next proposition identifies which parts of $\mathcal{T}_{\mathcal{A}}$ are associated with barrier curves that evolve backwards into $G_-$. To now numerically construct the admissible set one only needs to focus on these tangent points.

\begin{proposition}\label{prop:curves_A_in_G_minus_SEIR}
	Consider the constrained perfect SEIR model, \eqref{eq_SEIR_1}-\eqref{eq_SEIR_4}. The barrier $\DAM$ is constituted by the integral curves $x^{\bar{u}}$  satisfying \eqref{adjoint_SEIR_perfect}-\eqref{mini-beta-gamma_SEIR:eq}, that end tangentially to $G_{0}$ at a point of $\mathcal{T}_{\A}$, defined by \eqref{eq_tan_pts_SEIR}, for which 
	\begin{equation}\label{curves_A_in_G_minus_SEIR:eq}
	0\leq z_1\leq\min\{\frac{\gamma_{\max}}{\beta_{\min}}; 1 - \left( \frac{\gamma_{\max}}{\eta} + 1\right) I_{\max}\}
	\end{equation}
and is contained in $G_-$ at least for all $t\in]\bar{t} - \ee,\bar{t}[$, $\ee>0$.
\end{proposition}
\begin{proof}
	As in Proposition~\ref{prop:backwards_evolve}, the result follows from considering the inequality $\partial_- L_fg(x^{\bar{u}}(\bar{t}),\bar{u}(\bar{t}))< 0$ for points on $\mathcal{T}_{\A}$. First, evaluating $\dot{\lambda}_2(t) - \dot{\lambda}_1(t)$ at $t=\bar{t}$, thanks to \eqref{eq_tan_pts_SEIR} and \eqref{adjoint_SEIR_perfect}, we get $\dot{\lambda}_{2}(\bar{t}) -\dot{\lambda}_{1}(\bar{t})= -\eta < 0$. Therefore, $\lambda_2(t) - \lambda_1(t)$ is decreasing before vanishing at $\bar{t}$, which proves that $\lambda_2(t) - \lambda_1(t) > 0$ in an interval $] \bar{t}-\ee, \bar{t}[$, with $\ee>0$ suitably chosen. We immediately deduce that $\beta(t)= \beta_{\min}$ in the same time interval $] \bar{t}-\ee, \bar{t}[$. We also get, using the fact that $\lambda_{3}(\bar{t})=1$, that $\gamma = \gamma_{\max}$ in this interval (possibly with a smaller $\ee$). Thus $z_{2}(\bar{t})= \frac{\gamma_{\max}}{\eta}I_{\max}$. Moreover, in this interval, we must have
	$$\begin{aligned}
	\partial_- L_fg(x^{\bar{u}}(\bar{t}),\bar{u}(\bar{t}))& = \lim_{t\rightarrow \bar{t}, \,t\leq \bar{t}} \frac{d}{dt} \left (\eta E(t) - \gamma(t) I(t)\right) \\
	&= \eta I_{\max} \left( \beta_{\min} z_1 - \gamma_{\max}\right) < 0
	\end{aligned}$$
	which proves that $z_1$ must satisfy $z_1\leq \frac{\gamma_{\max}}{\beta_{\min}}$ and, thanks to the definition of $\mathcal{T}_{\A}$, $0\leq z_1\leq 1 - z_2 - z_3 = 1- \left( \frac{\gamma_{\max}}{\eta} - 1\right) I_{\max}$, hence \eqref{curves_A_in_G_minus_SEIR:eq}. The same Taylor's expansion as in Proposition~\ref{prop:backwards_evolve}, implies that the barrier curves evolve backwards into $G_-$ provided that \eqref{curves_A_in_G_minus_SEIR:eq} holds.
\end{proof}

\begin{remark}\label{barrier_extension_SEIR_rem}
Similar to Remark~\ref{barrier_extension:rem}, the set $\{(S,E,I): S\in [0,1], E = 0, I = 0\}$ is a set of equilibrium points for the SEIR model and clearly a subset of $\DAM$. Moreover, any integral curve of the system initiating from the set $\{(S,E,I): S = 0, E = 0, I \in [0, I_{\max}]\}$ results in $S(t) \equiv 0$, $E(t)\equiv 0$, $I(t) \leq I_{\max}$ for all $t\in[t_0,\infty[$, because $\gamma(t)>0$. Thus, as in the SIR case, the barrier $\DAM$ is the union of curves that satisfy the necessary conditions of Theorem~\ref{thm1} and this special subset that does not intersect $G_0$ tangentially.
\end{remark}

\section{Carrying out the Analysis of the MRPI for the Perfect SIR and SEIR Models (Problem~\ref{Pb2:def})}\label{section_MRPI_perfect_models}

In this section we carry out the detailed analysis of the Maximal Robust Positively Invariant sets of the constrained SIR and SEIR models under the assumption of perfect modelling (Problem~\ref{Pb2:def}).

As for the previous section, comparisons and numerical applications are presented in Sections~\ref{subsec_curves_backwards} and \ref{Sec_examples}.

\subsection{MRPI for Perfect SIR Model}\label{subsection_MRPI_SIR_perfect}

We consider the same setting as that of Subsection~\ref{subsection:perfect_model_case}, the constrained two-dimensional system \eqref{eq_SIR_1}-\eqref{eq_SIR_3}, and describe the system's MRPI set.

We first remark that the usable part $\DMO$, according to condition \eqref{MRPI-usable_part_eq}, is given by 
$$\max_{\beta \in [\beta_{\min}, \beta_{\max}]} L_f g = \beta SI - \gamma I \leq 0, \quad (S,I)\in G_{\Pi}~\mathrm{s.t.~} I-I_{\max}=0,$$
which yields $\DMO = \{(S,I) \mid 0\leq S \leq \min \left( \frac{\gamma}{\beta_{\max}}, 1-I_{\max}\right), \; I=I_{\max}\}$.
\bigskip

Now, carrying out the analysis with the ultimate tangentiality condition, \eqref{thm1_ult_tan_M}, and Hamiltonian maximisation condition \eqref{eq_Ham_max_M}, we note that the only difference with respect to \eqref{tancondSIRp:eq} is that we maximize with respect to $\beta$ in place of minimizing. Thus,
we find that the sole tangent point is located at $z=(\frac{\gamma}{\beta_{\max}}, I_{\max})$ with $0\leq \frac{\gamma}{\beta_{\max}} + I_{\max} \leq 1$. The adjoint equation being the same as \eqref{SIR_perfect_adjoint}, the input associated with the MRPI barrier curve is given by:
\begin{equation}
\bar{\beta}(t) = 
\begin{cases} 
\beta_{\max}&\mbox{if } \lambda_2(t) - \lambda_1(t) > 0, \\
\beta_{\min}&\mbox{if } \lambda_2(t) - \lambda_1(t) < 0, \\
\textrm{anything}&\mbox{if } \lambda_2(t) - \lambda_1(t) = 0,
\end{cases}\label{input_MRPI_SIR_perfect}
\end{equation}
which also only switches at isolated points in time, the proof being exactly the same as that of Proposition~\ref{prop_SIR_perfect_beta_non_singular}.

The following proposition is the analogue of Proposition~\ref{prop:backwards_evolve}. Its proof follows exactly the same lines and is left to the reader.
\begin{proposition}\label{prop:backwards_evolve_SIR_MRPI_perfect}
Consider the constrained perfect SIR model, \eqref{eq_SIR_1}-\eqref{eq_SIR_3}. The barrier $\DMM$ is constituted by the curve $x^{\bar{u}}$ generated by \eqref{input_MRPI_SIR_perfect} with \eqref{SIR_perfect_adjoint}, that ends tangentially to $G_0\cap \Pi$ at the point $(\frac{\gamma}{\beta_{\max}}, I_{\max})\in G_0\cap \Pi$ with $\frac{\gamma}{\beta_{\max}} + I_{\max} \leq 1$, and is contained in $G_-$ at least for all $t\in]\bar{t} - \ee,\bar{t}[$, $\ee >0$.
\end{proposition}
The reader may also verify that Remark~\ref{barrier_extension:rem} also applies to the barrier of the MRPI set. 
\subsection{MRPI for Perfect SEIR Model}\label{subsection_MRPI_perfect_SEIR}

We consider the same setting as that of Subsection~\ref{Subsection:admissible_set_SEIR}: the constrained three-dimensional system \eqref{eq_SEIR_1}-\eqref{eq_SEIR_4} and describe the system's MRPI set. The analysis follows exactly the same lines as in Subsection~\ref{Subsection:admissible_set_SEIR}, the operation of maximisation replacing the minimisation one. We therefore only sketch the results.

The usable part $\DMO$ is given by
$$\max_{\gamma\in [\gamma_{\min}, \gamma_{\max}]} L_f g = \eta E - \gamma I \leq 0, \quad (S,E,I)\in G_{\Pi}~\mathrm{s.t.~} I-I_{\max}=0.$$
Therefore,
$$\DMO = \{(S,E,I) \mid S  \in [0,1-I_{\max}], E \in [0, \min\{\frac{\gamma_{\min}}{\eta}I_{\max} , \; 1 - S - I_{\max}\}], I = I_{\max}\}.$$

\bigskip

Concerning the barrier, using conditions \eqref{eq_Ham_max_M} and \eqref{thm1_ult_tan_M}, we identify the line segment of potential tangent points:
\begin{equation}
\mathcal{T}_{\M_P}\triangleq\{ (z_1,z_2,z_3)^\top\in\mathbb{R}^3 : 0\leq z_1\leq 1 - z_2 - z_3; z_2 = \frac{\gamma_{\min}}{\eta} I_{\max}; z_3= I_{\max}\},\label{eq_tan_pts_M_SEIR}
\end{equation}
where the subscript $\M_{P}$ in $\mathcal{T}_{\M_P}$ refers to the ``perfect'' model case. 

The adjoint system being described by \eqref{adjoint_SEIR_perfect}, the input associated with the MPRI's barrier curves are thus given by:
\begin{equation}\label{maxi-beta-gamma_SEIR:eq}
\bar{\beta}(t) = 
\begin{cases} 
\beta_{\min}&\mbox{if } \lambda_2(t) - \lambda_1(t) < 0, \\
\beta_{\max}&\mbox{if } \lambda_2(t) - \lambda_1(t) > 0,\\
\textrm{anything}&\mbox{if } \lambda_2(t) - \lambda_1(t) = 0,
\end{cases}
~\mathrm{and~}
\bar{\gamma}(t) = 
\begin{cases} 
\gamma_{\min}&\mbox{if } \lambda_3(t) >  0, \\
\gamma_{\max}&\mbox{if } \lambda_3(t) < 0,\\
\textrm{anything}&\mbox{if } \lambda_3(t) = 0.
\end{cases}
\end{equation}
These inputs may only switch at isolated points in time provided that $z_{1}\neq 0$ in \eqref{eq_tan_pts_M_SEIR}, the proof following exactly the same arguments as the proof of Proposition~\ref{prop_SEIR_perfect_inputs_nonsingular}.  The next proposition is the analogue of Proposition~\ref{prop:curves_A_in_G_minus_SEIR}, its proof also following the same lines of reasoning.
\begin{proposition}\label{prop:curves_in_G_minus_SEIR_MPRI}
		Consider the constrained perfect SEIR model, \eqref{eq_SEIR_1}-\eqref{eq_SEIR_4}. The barrier $\DMM$ is constituted by the integral curves $x^{\bar{u}}$  satisfying \eqref{adjoint_SEIR_perfect} and \eqref{maxi-beta-gamma_SEIR:eq}, that end tangentially to $G_0$ at a point of $\mathcal{T}_{\M_P}$, defined by \eqref{eq_tan_pts_M_SEIR}, for which 
	\begin{equation}\label{curves_in_G_minus_SEIR:eq}
	0\leq z_1\leq\min\{\frac{\gamma_{\min}}{\beta_{\max}}; 1 - \left( \frac{\gamma_{\min}}{\eta} - 1\right) I_{\max}\}
	\end{equation}
and is contained in $G_-$ at least for all $t\in]\bar{t} - \ee,\bar{t}[$, $\ee>0$.
\end{proposition}
As before, the reader may verify that Remark~\ref{barrier_extension_SEIR_rem} applies to the MPRI for the perfect SEIR model as well.

\section{MRPI for the Imperfect SIR and SEIR Models (Problem~\ref{Pb3:def})}\label{section_MRPI_imperfect_models}

In this section we carry out the detailed analysis of the MRPIs of the constrained SIR and SEIR models under the assumption of imperfect modelling, corresponding to Problem~\ref{Pb3:def}).

\subsection{MRPI for Imperfect SIR Model}\label{subsubsec:MRPI_imperfect_SIR}

We again consider the normalised SIR model, but now assume that the model parameter $\gamma$ is \emph{not} perfectly known and that the intervention strategy, a simple affine feedback strategy, is given by:
\begin{equation}
\hat{\beta}(I(t)) = 
\beta_{\min}\left[\frac{I(t)}{I_{\max}}\right] + \beta_{\max}\left[1-\frac{I(t)}{I_{\max}}\right] = \frac{\beta_{\min}-\beta_{\max}}{I_{\max}} I(t) + \beta_{\max} ,\label{eq_beta_hat_feedback}
\end{equation}
where the level of social distancing depends on the proportion of infective individuals, has been pre-designed and implemented. Thus, the system is now:
\begin{align}
\dot S(t) &  =  - \hat{\beta}(t) S(t) I(t), \label{eq_SIR_uncertain_1}\\
\dot I(t) &  =  \hat{\beta}(t) S(t) I(t) - \gamma(t) I(t), \label{eq_SIR_uncertain_2}\\
I(t) - I_{\max} & \leq 0,\quad \gamma(t) \in [\gamma_{\min},\gamma_{\max}],\label{eq_SIR_uncertain_3}
\end{align}
with $t\in[t_0,\infty)$, and $0<\gamma_{\min}\leq\gamma_{\max}$. Again, to ease our notation, we will sometimes label the state $x\triangleq (S,I)^\top \in \mathbb{R}^2$, the uncertain model parameter $u\triangleq\gamma\in\mathbb{R}$, the parameter set $U \triangleq [\gamma_{\min}, \gamma_{\max}]\subset \mathbb{R}$, and $g(x) \triangleq I - I_{\max}$.

The adjoint equation \eqref{eq_adjoint} reads:
\begin{equation}\label{adjoint_SIR_imperfect}
\dot{\lambda} = 
\left(
\begin{array}{ccc}
\hat{\beta}(I)I & -\hat{\beta}(I)I \\
\alpha(I)S& -\alpha(I)S + \gamma
\end{array}
\right)\lambda,\,\,\lambda(\bar{t}) = (0,1)^T,
\end{equation}
where $\alpha(I)\triangleq 2\left(\frac{\beta_{\min} - \beta_{\max}}{I_{\max}}\right) I + \beta_{\max}$. Invoking conditions \eqref{eq_Ham_max_M} we find the input associated with the MRPI barrier curve, given by:
\begin{equation}\label{max_gamma_SIR_imperfect}
\bar{\gamma}(t) = 
\begin{cases} 
\gamma_{\min}&\mbox{if } \lambda_2(t) > 0, \\
\gamma_{\max}&\mbox{if } \lambda_2(t) < 0,\\
\textrm{anything}&\mbox{if } \lambda_2(t) = 0.
\end{cases}
\end{equation}
We again have a result concerning singular barrier curves.
\begin{proposition}
	Consider the constrained imperfect SIR model, \eqref{eq_SIR_uncertain_1}-\eqref{eq_SIR_uncertain_3}. Along any barrier integral curve satisfying \eqref{adjoint_SIR_imperfect}-\eqref{max_gamma_SIR_imperfect} and such that $z_1\neq 0$ and $2\beta_{\min}<\beta_{\max}$, the input $\bar{\gamma}$, given by \eqref{max_gamma_SIR_imperfect}, only switches at isolated points in time.
\end{proposition}
\begin{proof}
	 As before, assume that $\lambda_{2}(t)=0$ in some open interval $\mathcal{I}$. Then, indeed $\dot{\lambda}_{2}= 0= \alpha(I) S \lambda_1$. The reader may confirm that if $2\beta_{\min} < \beta_{\max}$, then $\alpha(I) \neq 0$. Thus, if also $z_1\neq 0$, we see that if $\lambda_2(t) = 0$ over an interval $\mathcal{I}$, then $\lambda_1(t) = \lambda_2(t) = 0$ on $\mathcal{I}$, hence the contradiction.
\end{proof}
 

From condition \eqref{thm1_ult_tan_M} the sole tangent point is located at $z=(\frac{\gamma_{\min}}{\beta_{\min}}, I_{\max})$, with $\frac{\gamma_{\min}}{\beta_{\min}} + I_{\max} \leq 1$, which we derive from the fact that $\hat{\beta}(I_{\max}) = \beta_{\min}$.

Again, the barrier curve evolves backwards into $G_-$:
\begin{proposition}\label{prop:backwards_evolve_SIR_MRPI_imperfect}
Consider the constrained imperfect SIR model, \eqref{eq_SIR_uncertain_1}-\eqref{eq_SIR_uncertain_3}.
	The barrier $\DMM$ is constituted by the curve $x^{\bar{u}}$ generated by  \eqref{max_gamma_SIR_imperfect} with \eqref{adjoint_SIR_imperfect}, that ends tangentially to $G_0\cap \Pi$ at the  point $(\frac{\gamma_{\min}}{\beta_{\max}}, I_{\max})\in G_0\cap \Pi$, with $\frac{\gamma_{\min}}{\beta_{\min}} + I_{\max} \leq 1$, and is contained in $G_-$ at least for all $t\in]\bar{t} - \ee,\bar{t}[$, $\ee>0$.
\end{proposition}
\begin{proof}
The reader may easily verify that $L_{f}g = \dot{I}(t) = \left(\hat{\beta}(t) S(t)  - \gamma(t) \right) I(t)$ and              $\frac{d}{dt} L_{f}g= -\hat{\beta}^{2}(t) I^{2}(t)S(t) +   \left(\hat{\beta}(t) S(t)  - \gamma(t) \right)\left( \beta_{\max}S(t)-\gamma(t) \right) I(t) - \dot{\gamma}(t) I(t)$. Thus, since $S(\bar{t})= \frac{\gamma_{\min}}{\beta_{\min}}$ and $I(\bar{t})= I_{\max}$, and since $\gamma(t)$ is constant in a suitable interval before $\bar{t}$, \ie~ $\dot{\gamma}(t) = 0$ in this interval, we have $\lim_{t\rightarrow \bar{t}, \,t\leq \bar{t}} \left( \hat{\beta}(t) S(t)  - \gamma(t) \right) = 0$ and
$\lim_{t\rightarrow \bar{t}, \,t\leq \bar{t}} \frac{d}{dt} L_{f}g=  -\hat{\beta}^{2}(\bar{t}) I_{\max}^{2}\frac{\gamma_{\min}}{\beta_{\min}} < 0$. The conclusion follows by the same Taylor's expansion argument as in Proposition~\ref{prop:backwards_evolve}.
\end{proof}

\subsection{MRPI for Imperfect SEIR Model}\label{MRPI_imperfect_SEIR}

Similar to the previous case, we consider the system:
\begin{align}
\dot S(t) &  =  - \hat{\beta}(t) S(t) I(t), \label{eq_SEIR_imperfect_1}\\
\dot E(t) &  =  \hat{\beta}(t) S(t) I(t) - \eta(t) E(t), \label{eq_SEIR_imperfect_2}\\
\dot I(t) &  =  \eta(t) E(t) - \hat{\gamma}(t) I(t), \label{eq_SEIR_imperfect_3}\\
I(t) - I_{\max} & \leq 0,\quad \eta(t) \in [\eta_{\min},\eta_{\max}],\label{eq_SEIR_imperfect_4}
\end{align}
with $t\in[t_0,\infty[$, $0<\eta_{\min}\leq \eta_{\max}$. The inputs $\hat{\beta}:[t_0,\infty[\rightarrow[\beta_{\min},\beta_{\max}]$ and $\hat{\gamma}:[t_0,\infty[\rightarrow[\gamma_{\min},\gamma_{\max}]$ are assumed to be pre-designed by the simple affine feedback strategies:
\begin{equation}\label{eq_feedback_beta}
\hat{\beta}(I(t)) = 
\beta_{\min}\left[\frac{I(t)}{I_{\max}}\right] + \beta_{\max}\left[1-\frac{I(t)}{I_{\max}}\right],
\end{equation}
and
\begin{equation}\label{eq_feedback_gamma}
\hat{\gamma}(I(t)) = 
\gamma_{\min}\left[1-\frac{I(t)}{I_{\max}}\right] + \gamma_{\max}\left[\frac{I(t)}{I_{\max}}\right]
\end{equation}
which may be interpreted as the fact that we increase social distancing and increase the rate of quarantining if the infective number increases. 

For this setting we identify: $x\triangleq (S,E,I)^\top \in \mathbb{R}^3$, the uncertain modelling parameter $u\triangleq\eta\in\mathbb{R}$, the uncertain parameter set $U \triangleq [\eta_{\min}, \eta_{\max}]\subset \mathbb{R}$, $g(x) \triangleq I - I_{\max}$, and $G_{\Pi}\triangleq G \cap \Pi$.

As before, the reader may easily verify that the set of potential tangent points is:
\begin{equation}
\mathcal{T}_{\M_I}\triangleq\{ (z_1,z_2,z_3)^\top\in\mathbb{R}^3 : 0\leq z_1\leq 1 - z_2 - z_3; z_2 = \frac{\gamma_{\max}}{\eta_{\max}} I_{\max}; z_3= I_{\max}\},\label{eq_tan_pts_M_SEIR_imperfect}
\end{equation}
where the subscript $M_{I}$ in $\mathcal{T}_{\M_I}$ denotes the ``imperfect'' case. The adjoint, \eqref{eq_adjoint}, reads:
\begin{equation}\label{adjoint_SEIR_imperfect}
\dot{\lambda} = 
\left(
\begin{array}{ccc}
\hat{\beta}(I)I & -\hat{\beta}(I)I & 0 \\
0 & \eta & -\eta \\
\alpha(I) S& -\alpha(I)S & \delta(I)
\end{array}
\right)\lambda,\,\,\lambda(\bar{t}) = (0,0,1)^T,
\end{equation}
where $\delta(I)\triangleq 2\left(\frac{\gamma_{\max} - \gamma_{\min}}{I_{\max}}\right) I + \gamma_{\min}$, and $\alpha(I)$ is as in Subsection~\ref{subsubsec:MRPI_imperfect_SIR}. 

It can be verified that the input associated with barrier curves of the MRPI is given by:
\begin{equation}\label{max_eta_SEIR_imperfect}
\bar{\eta}(t) = 
\begin{cases} 
\eta_{\max}&\mbox{if } \lambda_3(t) - \lambda_2(t) > 0, \\
\eta_{\min}&\mbox{if } \lambda_3(t) - \lambda_2(t) < 0,\\
\textrm{anything}&\mbox{if } \lambda_3(t) - \lambda_2(t) = 0.
\end{cases}
\end{equation}
\begin{proposition}
	Consider the constrained imperfect SEIR model, \eqref{eq_SEIR_imperfect_1}-\eqref{eq_SEIR_imperfect_4}. Along any barrier integral curve satisfying \eqref{adjoint_SEIR_imperfect}-\eqref{max_eta_SEIR_imperfect} and such that $z_1\neq 0$ and the determinant $\left\vert \begin{array}{cc}\alpha&\delta\\\dot{\alpha}&\dot{\delta}\end{array}\right\vert \neq 0$ for all $I\in [0, I_{\max}]$, the input $\bar{\eta}$, given by \eqref{max_eta_SEIR_imperfect}, only switches at isolated points in time.
\end{proposition}
\begin{proof}
	Assume, as before, that $\lambda_3(t) - \lambda_2(t)=0$ in a given interval of time $\mathcal{I}$. We immediately deduce that $\dot{\lambda}_{2}=0$ and $\dot{\lambda}_{3} - \dot{\lambda}_{2}=\dot{\lambda}_{3} = \alpha S (\lambda_{1}-\lambda_{2}) + \delta \lambda_{3} = 0$. Differentiating once more this latter expression yields $\ddot{\lambda}_{3}= \dot{\alpha}S(\lambda_{1}-\lambda_{2}) + \dot{\delta}\lambda_{3}$. Thus, if the above determinant does not vanish for all $I\in [0, I_{\max}]$ and if $S\neq 0$ we arrive at $\lambda_{1}=\lambda_{2}= \lambda_{3} = 0$, the desired contradiction.	
\end{proof}

The proof of the next proposition follows the same lines as before and is left to the reader.
\begin{proposition}\label{prop:curves_in_G_minus_imperfect_SEIR_MPRI}
Consider the constrained imperfect SEIR model, \eqref{eq_SEIR_imperfect_1}-\eqref{eq_SEIR_imperfect_4}. 
The barrier $\DMM$ is constituted by the curves $x^{\bar{u}}$ generated by  \eqref{max_eta_SEIR_imperfect} with \eqref{adjoint_SEIR_imperfect},
 that end tangentially to $G_0\cap \Pi$  at a point of $\mathcal{T}_{\M_I}$, defined by \eqref{eq_tan_pts_M_SEIR_imperfect}, with 
	\begin{equation*}
		0\leq z_1\leq\min\{\frac{\gamma_{\max}}{\beta_{\min}}; 1 - \left( \frac{\gamma_{\max}}{\eta_{\max}} - 1\right) I_{\max}\},
	\end{equation*}
and is contained in $G_-$ at least for all $t\in]\bar{t} - \ee,\bar{t}[$, $\ee>0$.
\end{proposition}

%
%
%
%
%
%

\section{Further Results on the Relative Location of the Admissible and MRPI Sets}\label{subsec_curves_backwards}

The admissible and/or MRPI sets may be trivially equal to the entire constrained state space, $G_{\Pi}$, in which case \emph{any} intervention strategy will maintain the infection cap. We now present some results that summarise the relationship between certain inequalities of the system parameters and the relative location of the sets.

\subsection{SIR model}\label{subsecSIR_curves_backwards}

\begin{proposition}\label{prop:nature_of_sets_perfect_model}
	Consider the constrained perfect SIR model \eqref{eq_SIR_1}-\eqref{eq_SIR_3}. We have the following:
	\begin{itemize}
		\item $\M = \A = G_{\Pi}$ if and only if $\beta_{\max} \leq \frac{\gamma}{(1 - I_{\max})}$,
		\item $\M \subsetneq \A = G_{\Pi}$, if and only if $\beta_{\max} > \frac{\gamma}{(1 - I_{\max})}$ and $\beta_{\min} \leq \frac{\gamma}{(1 - I_{\max})}$,
		\item $\M \subsetneq \A \subsetneq G_{\Pi}$, if and only $\beta_{\min} > \frac{\gamma}{(1 - I_{\max})}$.
	\end{itemize}
	Consider the constrained imperfect SIR model \eqref{eq_SIR_uncertain_1}-\eqref{eq_SIR_uncertain_3}, under the feedback \eqref{eq_beta_hat_feedback}. We have the following:
	\begin{itemize}
		\item $\M = G_{\Pi}$ if and only if $\beta_{\min} \leq \frac{\gamma_{\min}}{(1 - I_{\max})}$,
		\item $\M \subsetneq G_{\Pi}$ if and only if $\beta_{\min} > \frac{\gamma_{\min}}{(1 - I_{\max})}$.
	\end{itemize}
\end{proposition}
\begin{proof} The constraint set $G_{\Pi}$ is the polyhedron:
$$G_{\Pi} = \{ (S,I)\in \mathbb{R}^2 \mid I \leq I_{\max}, \; S\geq 0, \; I \geq 0, \; S + I \leq 1 \}$$
and its boundary is the union of the four faces
$$\partial G_{\Pi} = \{ (S,I)\in G_{\Pi} \mid  I = I_{\max}~\mathrm{or} \; S = 0~\mathrm{or} \; I = 0~\mathrm{or} \; S + I = 1 \}.$$
Indeed, $\M \subset \A \subset G_{\Pi}$ for every choice of $\beta_{\min}, \beta_{\max}, I_{\max}$.

Assume that  
\begin{equation}\label{betamax:ineq}
\beta_{\max} \leq \frac{\gamma}{(1 - I_{\max})}.
\end{equation} 
We will have $\M = G_{\Pi}$ if the vector field $f(S,I,\beta, \gamma)= \left( \begin{array}{c} -\beta SI\\\beta SI - \gamma I \end{array}\right)$ points towards the interior of  $G_{\Pi}$ along $\partial G_{\Pi}$ for all values of $\beta_{\min}, \beta_{\max}, I_{\max}$ satisfying \eqref{betamax:ineq}, \ie 
\begin{align}
&\max_{\beta\in[\beta_{\min},\beta_{\max}]} (0,1) f_{\big| I=I_{\max}} = (\beta_{max}S-\gamma)I_{\max} \leq 0\;  \; \forall S\in [0, 1-I_{\max}], \label{ineq1}\\
\mathrm{and} \; &\max_{\beta\in[\beta_{\min},\beta_{\max}]}  (1,0)f_{\big| S=0} = 0 \geq 0, \label{ineq2}\\
\mathrm{and} \; &\max_{\beta\in[\beta_{\min},\beta_{\max}]}  (0,1) f_{\big| I=0} =0  \geq 0, \label{ineq3}\\
\mathrm{and} \; &\max_{\beta\in[\beta_{\min},\beta_{\max}]}  (1,1) f_{\big| S+I=1} = -\gamma I \leq  0 \;\; \forall I \geq 0 \; \mathrm{such~that~} S+I=1. \label{ineq4}
\end{align}
The last three inequalities \eqref{ineq2}-\eqref{ineq4} are trivially satisfied and \eqref{ineq1} results from \eqref{betamax:ineq}, which proves that \eqref{betamax:ineq} implies $\M=G_{\Pi}$ and thus also $\M = \A = G_{\Pi}$. The converse, namely that $\M=G_{\Pi}$ implies \eqref{betamax:ineq}, follows from the fact that \eqref{ineq1} is valid only if \eqref{betamax:ineq} holds.

We proceed in the same way to prove that $\A = G_{\Pi}$ if $\beta_{\max} > \frac{\gamma}{(1 - I_{\max})}$ and $\beta_{\min} \leq \frac{\gamma}{(1 - I_{\max})}$, up to the fact that, in \eqref{ineq1}-\eqref{ineq4}, the maximisation with respect to $\beta$ is replaced by  the minimisation with respect to $\beta$.

The third bullet point results from the fact that none of the two previous cases hold, which prevents from having an equality of $\M$ or $\A$ to $G_{\Pi}$. 

The last two bullet points are obtained as the first and third bullets by changing $f$ in $\hat{f}(S,I, \gamma)= \left( \begin{array}{c} -\hat{\beta} SI\\\hat{\beta} SI - \gamma I \end{array}\right)$ and maximising with respect to $\gamma\in [\gamma_{\min},\gamma_{\max}]$.
\end{proof}

Proposition~\ref{prop:nature_of_sets_perfect_model} provides easily checkable inequalities of the system parameters that may be used to determine whether the sets are trivial or not. For example, focusing on the perfect SIR model, if $\beta_{\max}\leq \frac{\gamma}{(1-I_{\max})}$ then $G_{\Pi}$ is invariant for any choice of input and one need not worry about ever breaching the infection cap, as long as the state initiates inside $G_{\Pi}$. However, this is not the case when $\beta_{\max}>\frac{\gamma}{(1-I_{\max})}$, and the barrier of the MRPI must be computed by integrating the coupled differential equations of the system and the adjoint ones. If $\beta_{\min}>\frac{\gamma}{(1-I_{\max})}$ then the MRPI \emph{and} the admissible set are proper subsets of $G_{\Pi}$. In this case there are parts of $G_{\Pi}$ from which the infection cap will definitely be breached and, starting from $\A$ or $\M$, we must compute the corresponding barrier to derive an adapted intervention strategy.

\subsection{SEIR model}

We now summarise the analogous results for the perfect and imperfect SEIR models.
\begin{proposition}\label{prop:nature_of_sets_SEIR_model}
	Consider the constrained perfect SEIR model \eqref{eq_SEIR_1}-\eqref{eq_SEIR_4}. We have the following:
	\begin{itemize}
		\item $\M = \A = G_{\Pi}$ if and only if $\eta(1-I_{\max}) - \gamma_{\min}I_{\max}\leq 0$,
		\item $\M \subsetneq \A = G_{\Pi}$, if and only if $\eta(1-I_{\max}) - \gamma_{\min}I_{\max}> 0$ and $\eta(1-I_{\max}) - \gamma_{\max}I_{\max}\leq 0$,
		\item $\M \subsetneq \A \subsetneq G_{\Pi}$, if and only if $\eta(1-I_{\max}) - \gamma_{\max}I_{\max}> 0$.
	\end{itemize}
	Consider the constrained imperfect SEIR model \eqref{eq_SEIR_imperfect_1}-\eqref{eq_SEIR_imperfect_4} under the feedback \eqref{eq_feedback_beta} and \eqref{eq_feedback_gamma}. We have the following:
	\begin{itemize}
		\item $\M = G_{\Pi}$ if and only if $\eta_{\max}(1-I_{\max}) - \gamma_{\max}I_{\max}\leq 0$,
		\item $\M \subsetneq G_{\Pi}$ if and only if $\eta_{\max}(1-I_{\max}) - \gamma_{\max}I_{\max}> 0$.
	\end{itemize}
\end{proposition}
\begin{proof}
We follow the same arguments as in the proof of Proposition~\ref{prop:nature_of_sets_perfect_model}. For the SEIR model \eqref{eq_SEIR_1}-\eqref{eq_SEIR_4}, the set $G_{\Pi}$ is the polyhedron
$$G_{\Pi} = \{ (S,E,I)\in \mathbb{R}^3 \mid I \leq I_{\max}, \; S\geq 0, \; E\geq 0,\; I \geq 0, \; S +E+ I \leq 1 \}$$
and 
$$\partial G_{\Pi} = \{ (S,E,I)\in G_{\Pi} \mid  I = I_{\max}~\mathrm{or} \; S = 0~\mathrm{or} \; E = 0~\mathrm{or} \; I = 0~\mathrm{or} \; S +E+ I = 1 \}.$$
We thus prove that $\M=G_{\Pi}$ is equivalent to the first bullet point inequality, \ie
\begin{equation}\label{gammamax:ineq}
\eta (1-I_{\max}) - \gamma_{\min}I_{\max}\leq 0
\end{equation}
by showing that the vector field $f(S,E,I, \beta,\gamma) = \left( \begin{array}{c} -\beta SI\\ \beta SI - \eta E\\ \eta E -\gamma I \end{array}\right)$ (recall that $\eta$ is fixed) points towards the interior of $G_{\Pi}$ along $\partial G_{\Pi}$ for all values of $\beta\in [\beta_{\min}, \beta_{\max}]$, $\gamma\in [\gamma_{\min},\gamma_{\max}]$ provided that $ I_{\max}$ satisfies \eqref{gammamax:ineq}.
We thus evaluate the scalar product of $f$ with the normal to $\partial G_{\Pi}$ at each of its edges.
\begin{align*}
&\max_{\begin{array}{c}\beta\in[\beta_{\min},\beta_{\max}]\\ \gamma\in [\gamma_{\min},\gamma_{\max}]\end{array}} \max_{\begin{array}{c}E\in [0, 1-S-I_{\max}]\\ S\in [0,1]\end{array}}  (0,0,1) f_{\big| I=I_{\max}} \\
& = \max_{\gamma\in [\gamma_{\min},\gamma_{\max}]} \max_{\begin{array}{c}E\in [0, 1-S-I_{\max}]\\ S\in [0,1]\end{array}}  (\eta E -\gamma)I_{\max} \\
& = \eta (1- I_{\max}) - \gamma_{\min}I_{\max} \leq 0,
\end{align*}
and
\begin{align*}
&\max_{\begin{array}{c}\beta\in[\beta_{\min},\beta_{\max}]\\ \gamma\in [\gamma_{\min},\gamma_{\max}]\end{array}} \max_{\begin{array}{c}E\in [0, 1-I]\\ I\in [0,I_{\max}]\end{array}}  (1,0,0)f_{\big| S=0} =0 \geq 0,
\end{align*}
and
\begin{align*}
&\max_{\begin{array}{c}\beta\in[\beta_{\min},\beta_{\max}]\\ \gamma\in [\gamma_{\min},\gamma_{\max}]\end{array}} \max_{\begin{array}{c}S\in [0, 1-I]\\ I\in [0,I_{\max}]\end{array}}  (0,1,0) f_{\big| E=0} \\
&=  \max_{I\in [0,I_{\max}]} \beta I(1-I) = \begin{cases}
I_{\max}(1-I_{\max})&\mbox{if } I_{max}\leq \frac{1}{2}\\
\frac{\beta}{4}&\mbox{if } I_{max}\geq \frac{1}{2}
\end{cases} \\
& \geq 0
\end{align*}
and
\begin{align*}
&\max_{\begin{array}{c}\beta\in[\beta_{\min},\beta_{\max}]\\ \gamma\in [\gamma_{\min},\gamma_{\max}]\end{array}} \max_{\begin{array}{c}E\in [0, 1-S]\\ S\in [0,1]\end{array}}   (0,0,1) f_{\big| I=0} \\
&= \max_{\begin{array}{c}E\in [0, 1-S]\\ S\in [0,1]\end{array}} \eta E= \max_{S\in [0,1]} \eta (1-S) = \eta
\geq 0,
\end{align*}
and
\begin{align*}
&\max_{\begin{array}{c}\beta\in[\beta_{\min},\beta_{\max}]\\ \gamma\in [\gamma_{\min},\gamma_{\max}]\end{array}} \max_{\begin{array}{c}S+E\in [0, 1-I]\\ I\in [0,I_{\max}]\end{array}} (1,1,1) f_{\big| S+E+I=1} \\
&= \max_{I\in [0,I_{\max}]} -\gamma I = 0 \leq  0.
\end{align*}
Therefore, $f$ points towards the interior of $G_{\Pi}$ along $\partial G_{\Pi}$ for all values of $\beta$ and $\gamma$, or equivalently $\M=\A=G_{\Pi}$, if and only if \eqref{gammamax:ineq} holds.

We proceed in the same way to prove that $\A = G_{\Pi}$ if and only if $\eta(1-I_{\max}) - \gamma_{\min}I_{\max}> 0$ and $\eta(1-I_{\max}) - \gamma_{\max}I_{\max}\leq 0$, up to the fact that the maximisation with respect to $\beta$ and $\gamma$ is replaced by  the minimisation with respect to $\beta$ and $\gamma$.

As in the proof of Proposition~\ref{prop:nature_of_sets_perfect_model}, the third bullet point results from the fact that none of the two previous cases hold, which prevents from having an equality of $\M$ or $\A$ to $G_{\Pi}$. 

Finally, the last two bullet points are obtained as the first and third bullets by changing $f$ in $\hat{f}(S,E,I, \eta) = \left( \begin{array}{c} -\hat{\beta} SI\\ \hat{\beta} SI - \eta E\\ \eta E -\hat{\gamma} I \end{array}\right)$ and maximising with respect to $\eta\in [\eta_{\min},\eta_{\max}]$.
\end{proof}

\section{Examples}\label{Sec_examples}

We now present detailed numerical examples to clarify the results of the paper.

\subsection{SIR Examples}\label{subsection_SIR_examples}

\subsubsection{Perfect SIR Example}\label{subsubsection_SIR_perfect_Example}

First, consider the perfect SIR system, \eqref{eq_SIR_1}-\eqref{eq_SIR_3}, with $\gamma = 0.5$, $\beta_{\min} = 0.6$, $\beta_{\max} = 0.8$ and $I_{\max}\in\{0.02, 0.15, 0.4\}$. Thus, this (hypothetical) disease has an average recovery rate of 2 days and, on average, an individual nominally comes into contact (sufficiently to catch the disease) with $0.8$ other people per day. We have chosen these parameters to demonstrate the different types of possible sets. 

\begin{figure}[h]
	\begin{center}\
		\includegraphics[width=0.5\columnwidth]{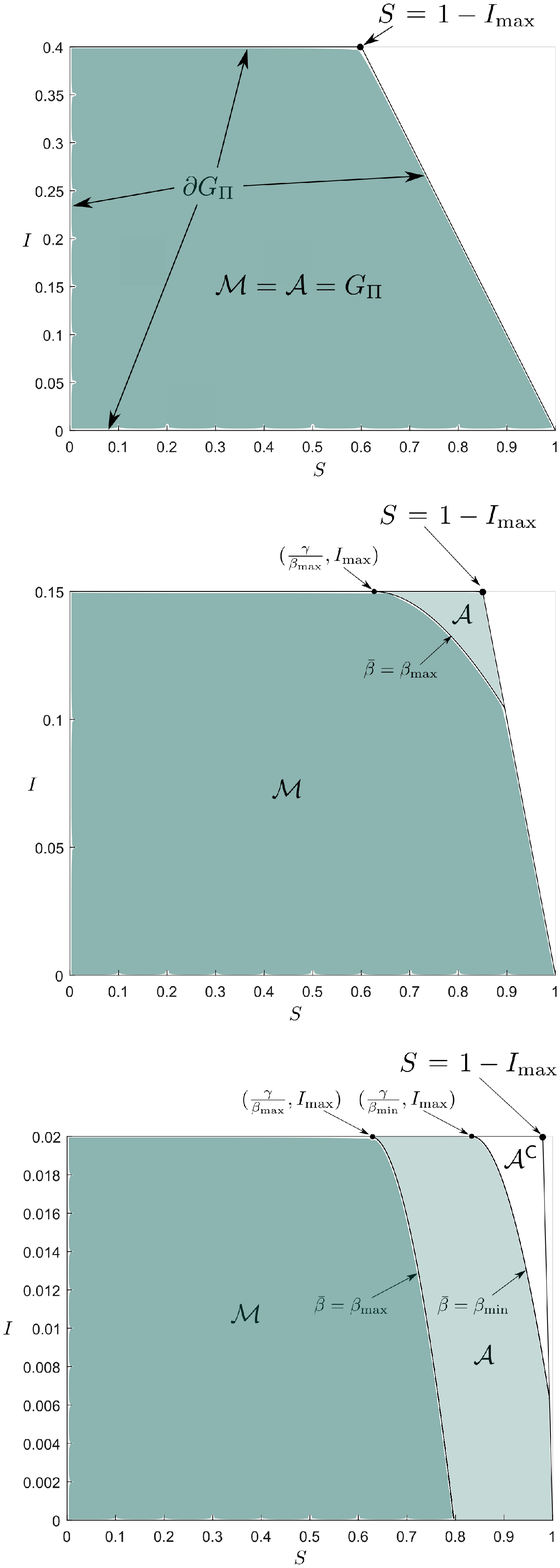} 
		\caption{\textbf{Sets for perfect SIR model example.} Projections of the admissible and MRPI sets onto the $S-I$ axes are shown. If the infection cap is large enough the entire constrained state space $G_{\Pi}$ is robustly invariant (top figure), otherwise the sets are subsets of $G_{\Pi}$.}
		\label{fig:SIR-perfect-model}
	\end{center}                                
\end{figure}

For $I_{\max} = 0.4$ we see from Proposition~\ref{prop:nature_of_sets_perfect_model} that $\M = \A = G_{\Pi}$. Thus, with such a large cap on the infection numbers the entire space $G_{\Pi}$ is robustly invariant. For $I_{\max}=0.15$ we see from Proposition~\ref{prop:nature_of_sets_perfect_model} that $\M \subsetneq \A = G_{\Pi}$. To now obtain the barrier associated with the MRPI, we integrate the system backwards in time from the tangent point $z=(\frac{\gamma}{\beta_{\max}},I_{\max})$ utilising the input \eqref{input_MRPI_SIR_perfect}. We stop integrating once the trajectory intersects the boundary of $G_{\Pi}$ because the set $\Pi$ is positively invariant (recall the discussions in Subsections~\ref{sec:SIR} and \ref{sec:SEIR}.) For $I_{\max} = 0.02$ we see from Proposition~\ref{prop:nature_of_sets_perfect_model} that $\M \subsetneq \A \subsetneq G_{\Pi}$. To obtain the barriers of the admissible set and MRPI we integrate the system backwards in time from the tangent points $z=(\frac{\gamma}{\beta_{\max}},I_{\max})$ and $z=(\frac{\gamma}{\beta_{\min}},I_{\max})$ utilising the input \eqref{input_MRPI_SIR_perfect} and \eqref{input_SIR_perfect_admissible}, respectively. See Figure~\ref{fig:SIR-perfect-model} for details.


It turns out that the input associated with the admissible set is saturated at $\beta_{\min}$, and the one associated with the MRPI is saturated at $\beta_{\max}$. Moreover, we can see that, going backwards, the barrier curves always evolve into $G_-$, as is expected from Propositions~\ref{prop:backwards_evolve} and \ref{prop:backwards_evolve_SIR_MRPI_perfect}.

\begin{figure}[h]
	\begin{center}\	\includegraphics[width=0.5\columnwidth]{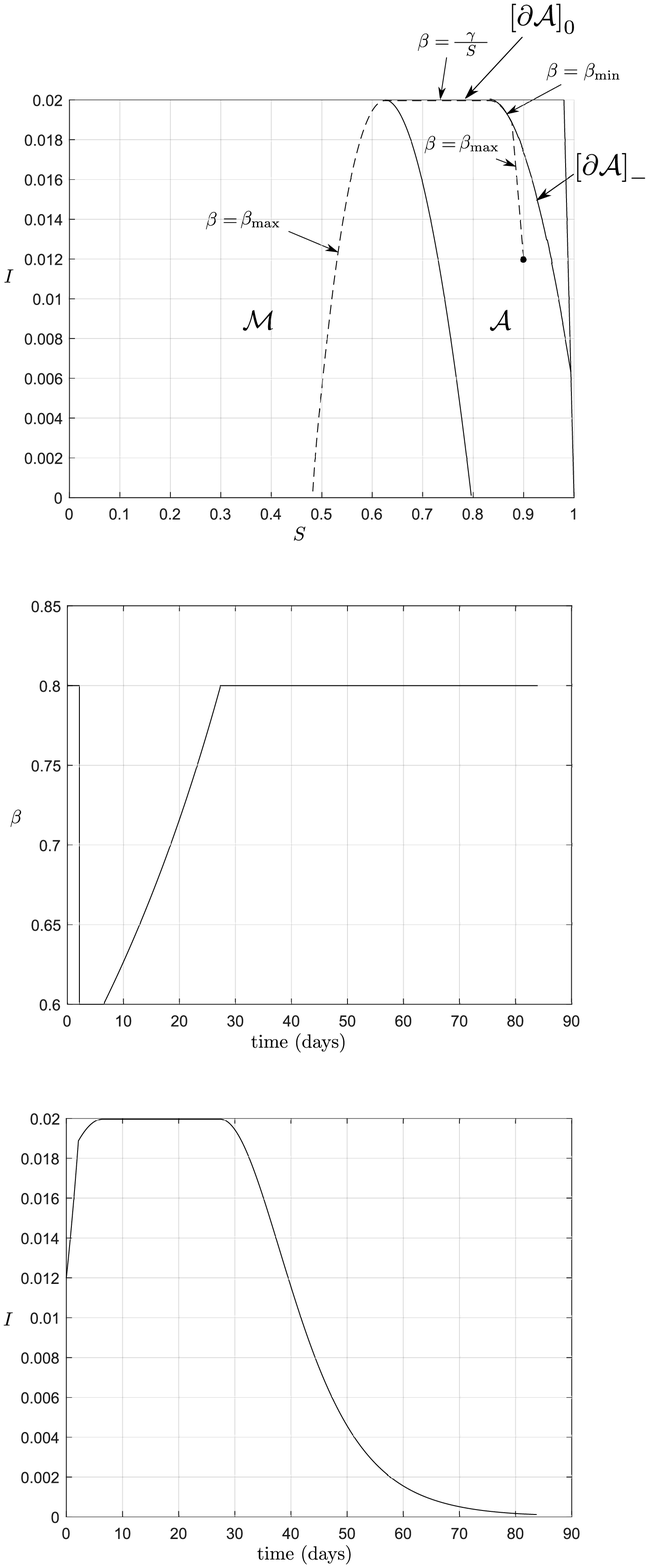} 
		\caption{\textbf{Solution with the switching law.} The state trajectory is shown in the top figure (dashed curve). The middle plot is of the intervention strategy, $\beta$, and the right plot is the proportion of infectives over time.}
		\label{fig:SIR-perfect-model_2}                             
	\end{center}                                
\end{figure}

To now use the sets in the management of an epidemic the intervention strategy may be chosen according to the location of the state. A possible strategy could be:
\begin{itemize}
	\item If $(S(t), I(t))^\top \in \M \cup \mathsf{int}(\A)$, let $\beta(t) = \beta_{\max}$,
	\item If $(S(t), I(t))^\top \in \partial \A$, let $\beta(t) = \beta_{\min}$,
	\item If $(S(t), I(t))^\top \in \DAO$, let $\beta(t)= \frac{\gamma}{S(t)}$.
\end{itemize}
To clarify, if the state is located in the MRPI then it is guaranteed that the infection cap can always be maintained, and the nominal contact rate, $\beta_{\max}$, may be allowed. The same is true if the state is located in the interior of the admissible set, but this may only be possible for some period of time. If the state reaches the admissible set's barrier, $\DAM$, then maximal social distancing must be imposed, thus $\beta(t) = \beta_{\min}$, and the infection cap will be reached in finite time, but never breached. If the state reaches the admissible set's usable part, $\DAO$, then some freedom in the social distancing is allowed. In fact, we may choose $\beta(t)$ such that $\dot{I} = 0$, which yields $\beta(t) = \frac{\gamma}{S(t)}$. Lastly, the set $\Acomp$ should always be avoided, as from here the infection cap is guaranteed to be breached regardless of $\beta$. Figure~\ref{fig:SIR-perfect-model_2} shows the result of this switching law being applied from an arbitrary initial condition, $(S(t_0),I(t_0))^{\top} = (0.8,0.012)^{\top}\in\mathsf{int}(\A)$, for the case where $I_{\max} = 0.02$.

\subsubsection{Imperfect SIR Example}

\begin{figure}[h]
	\begin{center}\	\includegraphics[width=0.5\columnwidth]{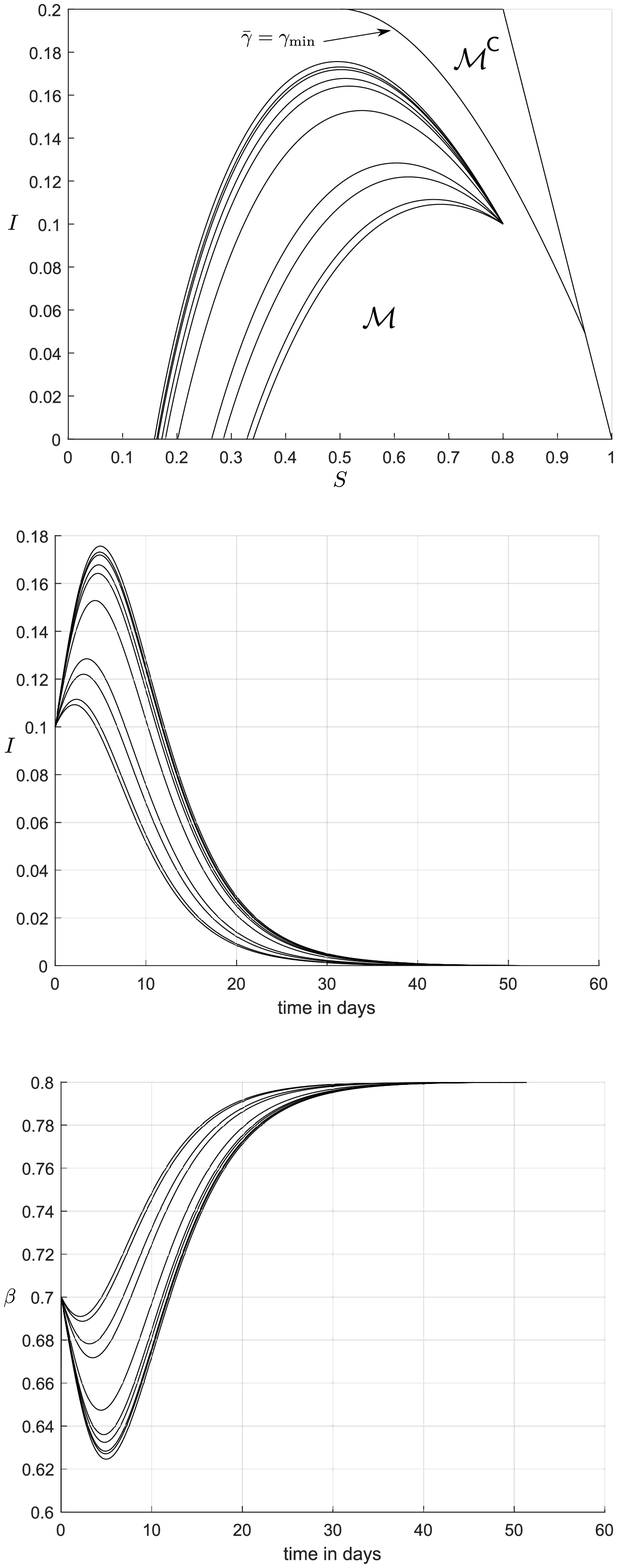} 
		\caption{\textbf{MRPI for imperfect SIR model example.} We also show ten simulations from an initial state $x_0\in\M$ with randomly sampled $\gamma\in[0.3,0.5]$.}
		\label{fig:SIR-imperfect-model}                             
	\end{center}                                
\end{figure}

Consider the imperfect SIR model, \eqref{eq_SIR_uncertain_1}-\eqref{eq_SIR_uncertain_3}, with the simple feedback strategy $\hat{\beta}$ as in \eqref{eq_beta_hat_feedback} and the constants $I_{\max} = 0.2$, $\beta_{\min} = 0.6$, $\beta_{\max} = 0.8$, and $\gamma_{\min} = 0.3$, $\gamma_{\max} = 0.5$. Thus, the recovery rate for this hypothetical disease may vary. We see from Proposition~\ref{prop:nature_of_sets_perfect_model} that $\mathcal{M}\subsetneq G_{\Pi}$, and so we find the barrier of the system's MRPI according to the analysis in Subsection~\ref{subsubsec:MRPI_imperfect_SIR}, see Figure~\ref{fig:SIR-imperfect-model}. It turns out that $\bar{\gamma}$ is saturated at $\gamma_{\min}$ all along the MRPI's barrier. We also show a few solutions of the system from the initial state $(S(t_0),I(t_0))^{\top} = (0.8,0.1)^\top\in\M$ with $\gamma\in [0.3,0.5]$, emphasising that the infection cap can be maintained with this feedback regardless of the modelling uncertainty.

\begin{remark}
	In epidemiology it is common to want an intervention strategy that renders the \emph{basic reproduction number}, $\mathcal{R}_0$, less than one in order to eradicate the disease. However, this is not required to maintain the infection cap using the introduced sets. In fact, in the examples of Subsection~\ref{subsection_SIR_examples}, $\mathcal{R}_0 = \frac{\beta}{\gamma} > 1$ for all mentioned values of $\beta$ and $\gamma$.
\end{remark}

\subsection{SEIR Examples}\label{subsec_SEIR_example}

\subsubsection{Perfect SEIR Example}

\begin{figure}[h]
	\begin{center}\	\includegraphics[width=\columnwidth]{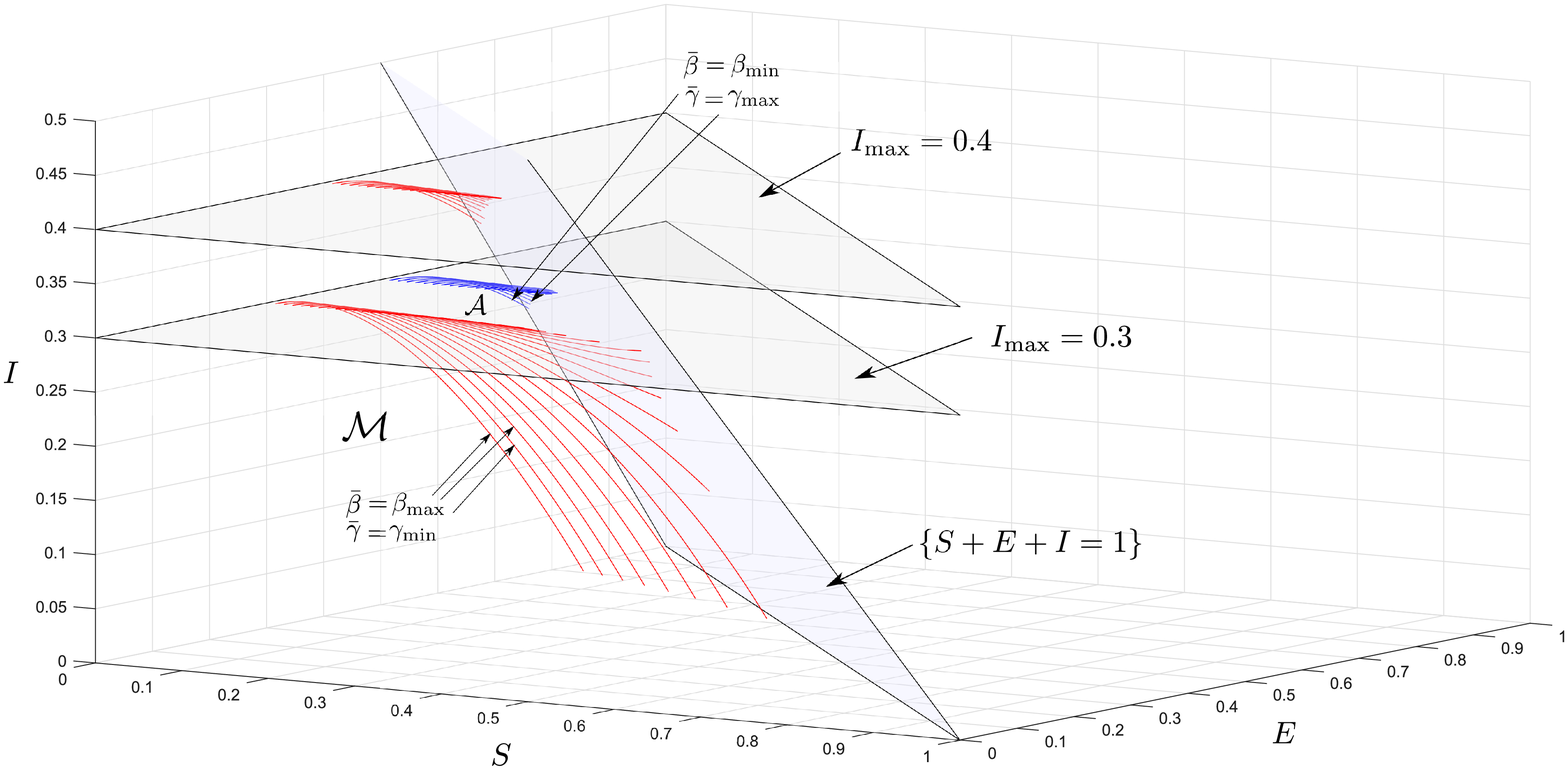} 
		\caption{\textbf{Sets for perfect SEIR model, with $I_{\max}\in\{0.3,0.4\}$.} Shown are curves running along the barriers of the MRPI (in red) and the admissible set (in blue).} 
		\label{fig:SEIR_perfect_03_04}                             
	\end{center}                                
\end{figure}

Consider the perfect SEIR model, \eqref{eq_SEIR_1}-\eqref{eq_SEIR_4}, with the constants $\beta_{\min} = 0.8$, $\beta_{\max} = 1$, $\gamma_{\min} = \frac{1}{5}$, $\gamma_{\max} = \frac{1}{3}$, $\eta = \frac{1}{5}$ and $I_{\max} \in\{0.3,0.4\}$. From Proposition~\ref{prop:nature_of_sets_SEIR_model} we see that $\mathcal{M}\subsetneq \mathcal{A} = G_{\Pi}$ for $I_{\max} = 0.4$. Thus, we sample final tangent points along the line segment $\mathcal{T}_{\M_{P}}$ for which $z_1 \leq \min\{\frac{\gamma_{\min}}{\beta_{\max}}; 1 -  \gamma_{\min}I_{\max}/\eta - I_{\max}\}$ (see Proposition~\ref{prop:curves_in_G_minus_SEIR_MPRI}) and integrate backwards with the appropriate inputs, as specified in Subsection~\ref{subsection_MRPI_perfect_SEIR}, terminating the integration if the curves intersect the infection cap, or the plane $\{(S,E,I)\in\mathbb{R}^3; S+E+I = 1\}$. This produces the red curves in Figure~\ref{fig:SEIR_perfect_03_04}, running along the MRPI's barrier.

For $I_{\max}= 0.3$ we see that $\mathcal{M}\subsetneq \mathcal{A} \subsetneq G_{\Pi}$. Thus, for this case we find both $\mathcal{M}$ and $\mathcal{A}$, integrating backwards from $\mathcal{T}_{\M_{P}}$ and $\mathcal{T}_{\A}$ with the appropriate inputs as given in Subsections~\ref{subsection_MRPI_perfect_SEIR} and \ref{Subsection:admissible_set_SEIR}.

The inputs $\bar{\beta}$ and $\bar{\gamma}$ are saturated at $\beta_{\max}$ and $\gamma_{\min}$ (resp. $\beta_{\min}$ and $\gamma_{\max}$) for the barrier of the MRPI (resp. the barrier of the admissible set). As was observed in the examples for the SIR model, if $I_{\max}$ is large enough the entire set $G_\Pi$ is robustly invariant, and if $I_{\max}$ is small enough, the sets $\A$ and $\M$ become nontrivial.

\subsubsection{Imperfect SEIR Example}

Now consider the imperfect SEIR model, \eqref{eq_SEIR_imperfect_1}-\eqref{eq_SEIR_imperfect_4} with the parameters $\beta_{\min}=0.8$, $\beta_{\max} = 1$, $\gamma_{\min} = \frac{1}{5}$, $\gamma_{\max} = \frac{1}{3}$, $\eta_{\min} = \frac{1}{7}$, $\eta_{\max} = \frac{1}{5}$ and $I_{\max} =  0.1$, and the feedback \eqref{eq_feedback_beta} and \eqref{eq_feedback_gamma}. From Proposition~\ref{prop:nature_of_sets_SEIR_model}, we have $\mathcal{M}\subsetneq G_{\Pi}$. Thus, we integrate backwards from points in $\mathcal{T}_{\mathcal{M}_I}$ for which $z_1\leq \min\{\frac{\gamma_{\max}}{\beta_{\min}}; 1 -  \gamma_{\max}I_{\max}/\eta_{\max} - I_{\max}\}$ (see Proposition~\ref{prop:curves_in_G_minus_imperfect_SEIR_MPRI}) with the input $\bar{\eta}$ as specified in Subsection~\ref{MRPI_imperfect_SEIR}.

\begin{figure}[h]
	\begin{center}\
		\includegraphics[width=\columnwidth]{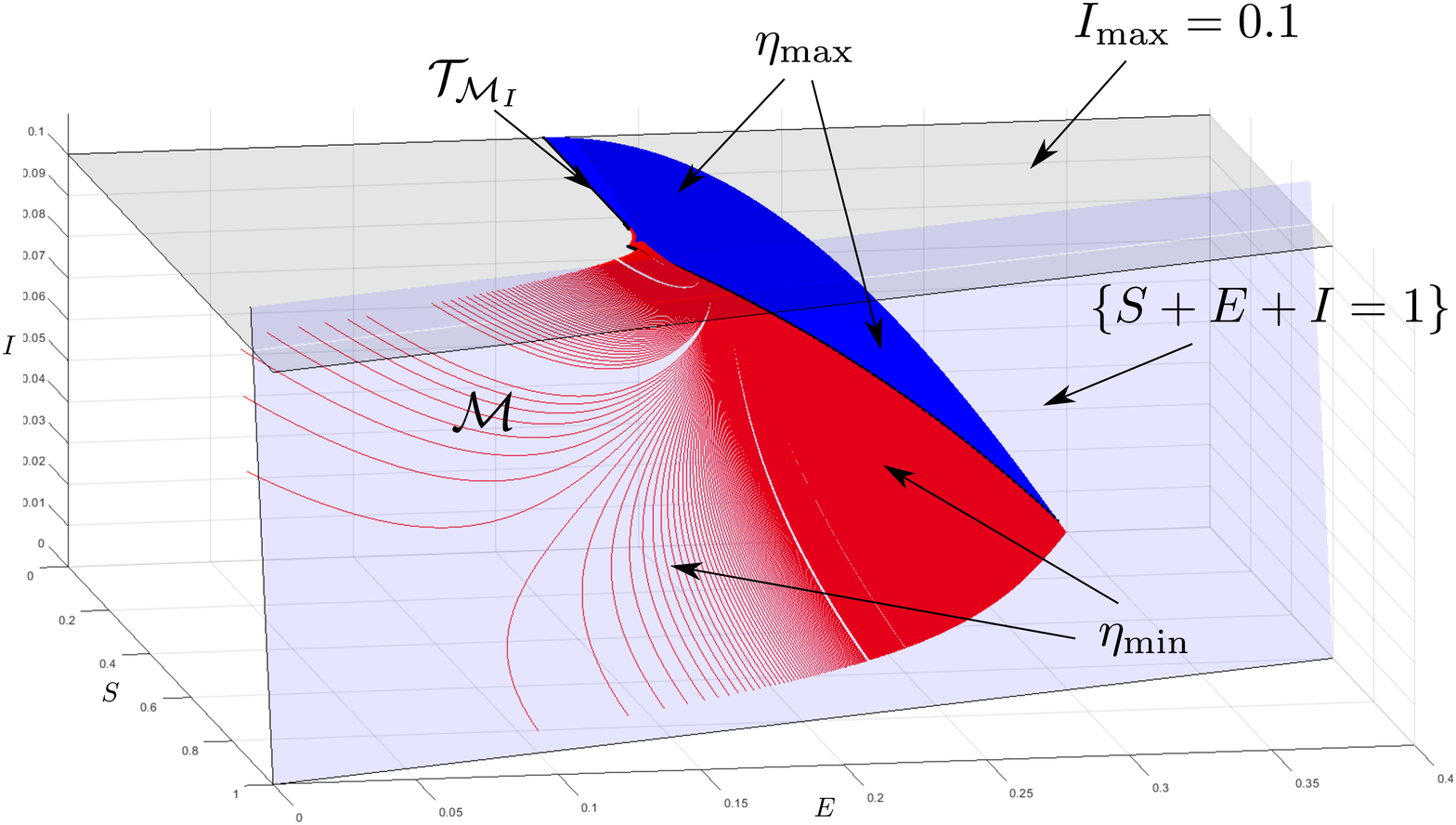} 
		\caption{\textbf{MRPI for imperfect SEIR model.} Going backwards in time, the input associated with the barrier, $\bar{\eta}$, switches from $\eta_{\max}$ (blue curves) to $\eta_{\min}$ (red curves) at the black dots.}
		\label{fig:SEIR_imperfect}                             
	\end{center}                                
\end{figure}

An unexpected result for this set is that the input associated with barrier curves, $\bar{\eta}$, switches from $\eta_{\min}$ to $\eta_{\max}$ at some point along the barrier (these points are indicated on the figure). To explain this, we can interpret $\eta$ as an input with the goal of maximising $I$ over time. For this example it is optimal for it to remain at $\eta_{\min}$ in order to build up the proportion of exposed individuals to some threshold when it switches to $\eta_{\max}$, which then results in a large increase in the proportion of infectives. This is the worst-case change in $\eta$ that can be expected. The computed set is shown in Figure~\ref{fig:SEIR_imperfect}.

\section{Discussion}\label{sec_perspectives}

We applied the theory of barriers to describe integral curves of the system that run along special parts of the boundaries of the admissible and MRPI sets of the constrained SIR and SEIR models. The analysis in Sections~\ref{SIR_analysis}, \ref{section_MRPI_perfect_models} and \ref{section_MRPI_imperfect_models} summarises the details required to construct the sets for each of the three problems stated in Subsection~\ref{Sec_use_of_sets}, with Section~\ref{subsec_curves_backwards} presenting inequalities of the system parameters that allow one to determine whether the sets are trivially equal to the constrained state space $G_{\Pi}$ or not. We also demonstrated, in the example in Subsection~\ref{subsection_SIR_examples} how an intervention strategy may be chosen using the computed sets in order to maintain the infection cap.

There are many avenues of future research to pursue: the analysis could be applied to other epidemic models; we could investigate the efficacy of combining model predictive control with knowledge of the sets; or look into the effects of more realistic intervention strategies that must remain constant over long periods of time, as opposed to the continuous ones introduced in \eqref{eq_beta_hat_feedback} or \eqref{eq_feedback_gamma} - \eqref{eq_feedback_beta}.

\section*{Acknowledgments}

This work was partially funded by the German Federal Ministry of Education and Research (BMBF; grants 05M18OCA: "Verbundprojekt 05M2018 - KONSENS: Konsistente Optimierung und Stabilisierung elektrischer Netzwerksysteme")

\nolinenumbers

\bibliographystyle{plos2015}
\bibliography{references}

\end{document}